\documentclass[11pt]{article}
\pdfoutput=1

\usepackage{lineno}
\usepackage{color}
\usepackage{amsfonts}
\usepackage{amsmath}
\usepackage{amssymb}
\usepackage{amsthm}
\usepackage{subeqnarray}
\usepackage{lscape,graphicx,url}
\usepackage{anysize}
\usepackage{proof}
\usepackage{fancyhdr}
\usepackage{latexsym}
\usepackage{array}
\usepackage{multirow}
\usepackage{textcomp}
\usepackage{optprog}
\usepackage{mathtools}
\usepackage{natbib}

\usepackage{enumitem}
\usepackage{csquotes}
\usepackage{hyperref}

\newtheorem{lemma}{Lemma}

\title{Optimally solving the joint order batching and picker routing problem}

\author{Cristiano Arbex Valle$^1$\thanks{Cristiano Arbex Valle is funded by CPNq grant 401367/2014-2.} \and John E Beasley$^2$ \and Alexandre Salles da Cunha$^1$\thanks{Alexandre Salles da Cunha is partially funded by FAPEMIG grant CEX - PPM-00187-15 and CNPq grants 303677/2015-5, 471464/2013-9 and 200493/2014-0.}}

\date{}

\begin{document}

\maketitle

\begin{center} 
{\footnotesize

$^1$Departamento de Ci\^{e}ncia da Computa\c{c}\~{a}o, \\
 Universidade Federal de Minas Gerais, \\
 Belo Horizonte, MG 31270-010, Brasil \\
\{arbex,acunha\}@dcc.ufmg.br \\ \vspace{0.3cm}
$^2$Brunel University \\ Mathematical Sciences, UK \\ john.beasley@brunel.ac.uk \\
}
\end{center}

\begin{abstract}
In this work we investigate the problem of order batching and picker routing in storage areas. These are labour and capital intensive problems, often responsible
for a substantial share of warehouse operating costs. In particular, we consider the case of online grocery shopping in which orders may be composed of dozens of items. 

We present a formulation for the problem based on an exponential number of connectivity constraints and we introduce a significant number of valid inequalities based on the standard layout of warehouses, composed of parallel aisles and two or more cross-aisles. The proposed inequalities are highly effective and greatly improve computational results. 

Instances involving up to 20 orders are solved to proven optimality when we jointly consider order batching and picker routing. Instances involving up to 5000 orders are considered where order batching is done heuristically, but picker routing is done optimally.
\end{abstract}


\section{Introduction}

Warehouses require intensive product handling operations, amongst those order picking is known to be the most labour and machine intensive. Its cost is estimated to be as much as 55\% of total warehouse operating expenses \citep{tompkins2010}; as online shopping has grown in popularity in recent years we believe that this figure is likely increasing. Order picking is defined as being \textit{the process of retrieving products from storage in response to specific customer requests}. 

Material handling activities can be differentiated as parts-to-picker systems, in which automated units deliver the items to stationary pickers, and picker-to-parts systems, in which pickers walk/ride through the warehouse collecting requested items. With respect to the latter, \cite{henn2012} distinguish three planning problems: assignment of products to locations, grouping of customer orders into batches, routing of order pickers. This paper deals jointly with the last two activities, which are crucial to the efficiency of warehouse operations.

In particular, we consider the case of online grocery shopping, where orders may be composed of dozens of items. As practical aspects of this problem we highlight (i) the heterogeneity of products (which can be of various shapes, sizes or expiration dates) and (ii) the fact that the picking can be performed in warehouses that are either closed or open to the public (such as supermarkets); due to such features pickers generally walk/ride the warehouse, collecting products manually instead of relying on automated systems. Pickers place products in baskets; when orders can be composed of dozens of items, mixing different orders in the same basket or splitting an order among different pickers is generally avoided to reduce order processing errors.

We formulate and solve the Joint Order Batching and Picker Routing Problem (JOBPRP). The task is to find minimum-cost closed walks 
 where each picker visits all locations required to pick all products from their assigned orders. Locations may be visited more than once if necessary. In our previous work \citep{valle2016} we introduced three integer programming formulations for JOBPRP, each one being the basis for a different branch-and-bound algorithm. One of the formulations is a directed model that involves exponentially many constraints to enforce connectivity requirements for closed walks. The other two are compact formulations based on network flows. We introduced a branch-and-cut algorithm that relies on the non-compact model and examined the compact formulations using the CPLEX branch-and-bound solver. Results suggested that the non-compact formulation and the compact formulation based on single-commodity network flows are competitive. In our experiments, the formulation based on multi-commodity flows ran out of memory for larger instances.

In this work, we focus on an improved version of the non-compact formulation. Our main contribution is the introduction of several valid inequalities (cuts) based on a sparse graph representation of warehouses which, to the best of our knowledge, have never been proposed before. We show that the inclusion of these cuts in the non-compact formulation greatly improves computational results. 
The proposed cuts prevent many solutions that violate subtour breaking constraints and thus the (relatively) expensive max-flow separation algorithm can be switched off. Moreover, if the batching and routing problems are solved separately, we show that, given an assignment of orders to pickers, an optimal routing can be computed very quickly. We also introduce a special case of the problem where pickers are only allowed to reverse at corners, which is more intuitive for humans, reduces computation times and still produces solutions of very good quality when compared to the original problem.

As in \cite{valle2016}, we also introduce a JOBPRP test instance generator partially based on publicly available real-world data. The test instances used in this work are available for 
download\footnote{Test instances available at \href{http://www.dcc.ufmg.br/~arbex/orderpicking.html}{http://www.dcc.ufmg.br/$\sim$arbex/orderpicking.html}\label{footnote1}}.


The remainder of this paper is organised as follows. Section~\ref{sec:lit} presents a description and a literature review of the problem. Section~\ref{sec:graph} introduces JOBPRP as a graph optimisation problem and Section~\ref{sec:IP} presents the 
aforementioned
non-compact integer programming formulation  and briefly discusses some implementation details. Section~\ref{sec:ineq} introduces several classes of valid inequalities based on the standard layout of warehouses. Section~\ref{sec:results} discusses computational results whilst Section~\ref{sec:application} discusses the application of the approach in this paper in the context of supermarket home deliveries.
Finally 
in Section~\ref{sec:conclusion} we present some concluding remarks and future research directions.

\section{Problem description and literature review}
\label{sec:lit}

\begin{figure}[!htb]
\centering
  \includegraphics[width=0.5\textwidth]{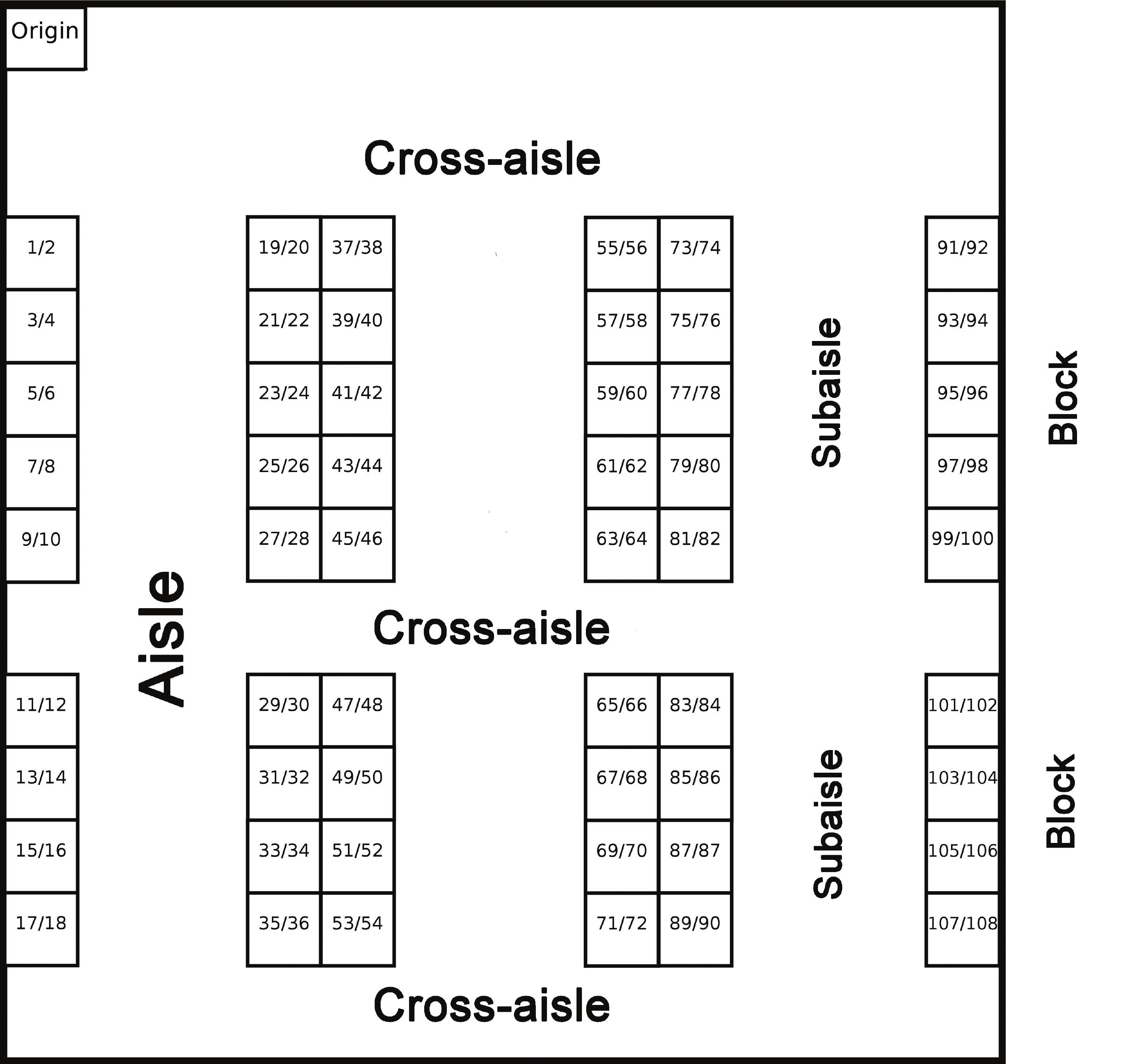}
  \caption{Warehouse structure}
  \label{fig1}
\end{figure}

A warehouse has a rectangular layout with no unused space and consists of parallel \textit{aisles}. The warehouse is divided into a number of \emph{blocks}. Blocks are separated by \emph{cross-aisles}. A \emph{subaisle} is defined as a section of an aisle within a block. Every warehouse has at least one cross-aisle at the top and one at the bottom, but may contain more. Figure \ref{fig1} illustrates these concepts. 

Aisles contain slots on both sides, each slot holds one type of product. Slots can also be stacked vertically in shelves. We assume that pickers move in the centre of an aisle and that products on both sides can be reached by the picker. 

\begin{figure}[!htb]
\centering
\begin{minipage}[b]{0.49\textwidth}
    \includegraphics[width=1.18\textwidth]{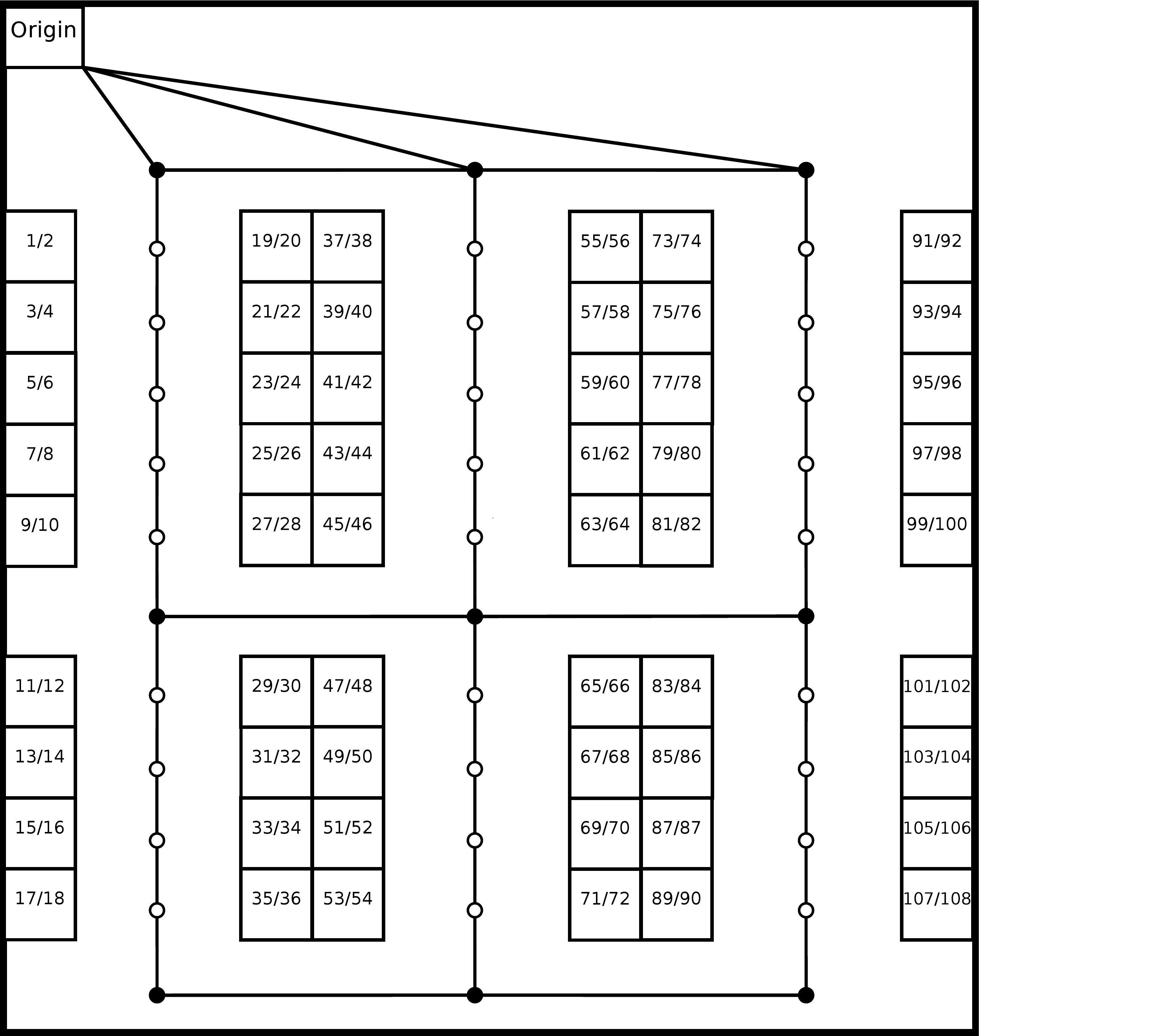}
    \caption{Example of warehouse/graph layout}
    \label{fig2}
\end{minipage}
\hfill
\begin{minipage}[b]{0.49\textwidth}
    \includegraphics[width=1.18\textwidth]{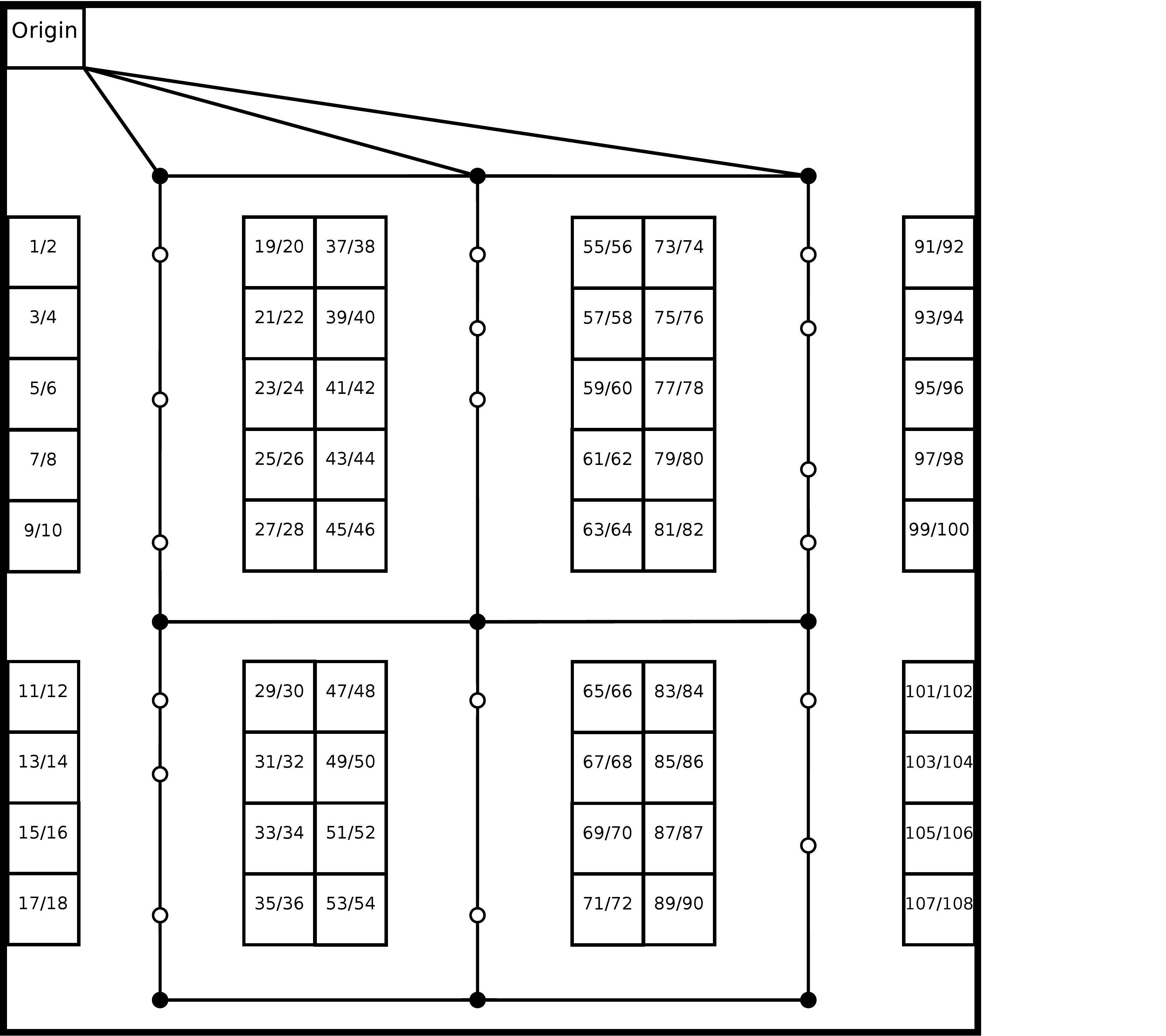}
    \caption{Warehouse with reduced graph}
    \label{fig3}
\end{minipage}
\end{figure}

A typical warehouse layout with three cross-aisles, three aisles  and two shelves, together with a sparse graph representation, can be seen in Figure \ref{fig2}. 
White vertices represent locations from where the picker reaches products, i.e. from the top left white vertex the picker reaches products in slots 1, 2, 19 and 20. 
Note here that whilst in Figure \ref{fig2} we have shown just two slots on each side of an aisle this is purely for ease of illustration. Our approach is not dependent on the number of such slots and can be applied without modification to problems with more slots.

In Figure \ref{fig2} there are six subaisles, each consisting of four consecutive white vertices. 
Black vertices represent ``artificial'' locations which connect aisles and cross-aisles, while the ``Origin'' indicates the starting and return point for pickers. Often, however, not all locations need to be visited as the joint set of orders does not contain products in all locations. It is therefore possible to reduce the graph that depicts the warehouse by eliminating such locations, as it can be seen in Figure \ref{fig3}, where, for example, no products in slots 3, 4, 21 and 22 are needed.

Pickers use trolleys that accommodate a limited number of baskets\footnote{In different industries, the terms cart and tote are more often used than trolley and basket.
}. The number of baskets necessary to carry each order is assumed to be known. In fact, in this work, we consider that all baskets carry a fixed number of items, irrespective of their shapes and sizes. Baskets needed to carry a batch must not exceed trolley capacity. Following general practice in supermarkets, orders from different customers cannot be put together in the same basket and a single order cannot be split between different trolleys. Mixing and dividing orders could reduce picking time, however this benefit is often offset by necessary post-processing and a higher risk of errors.

There is a vast literature on the problems of order batching and picking. An extensive survey relating to different picking strategies was presented by \cite{deKoster2007}, where order picking is shown to be the most critical activity for productivity growth in the sector. Wave picking is an example of an alternative manual picking strategy, where the warehouse is divided into zones, each zone being assigned to one or more pickers; however in this strategy post-processing of orders is necessary. In this paper we deal with batch picking, where a picker collects all of the products for one or more orders.

The special case where a single picker is sufficient to collect all orders is a particular case of the Travelling Salesman Problem (TSP). \cite{ratliff1983} showed that for warehouses which have only two cross-aisles (so a single block), the problem can be solved polynomially via a dynamic programming approach. The algorithm by Ratliff \& Rosenthal is based on a sparse graph representation of the problem, whereas previous works on the TSP generally assumed a complete graph. The sparse graph representing a warehouse containing a single block is a type of \textit{series parallel graph}. \cite{cornuejols1985} generalised the algorithm by Ratliff \& Rosenthal and showed that the TSP in any series parallel graph is solvable in polynomial time. In that work, the authors coined the term Graphical TSP (GTSP) to characterise the TSP problem on sparse graphs - a major difference being that arcs and vertices can be revisited in a tour, unlike the TSP in a complete graph. \cite{fonlupt1992} showed that the GTSP is polynomially solvable whenever the corresponding graph is {\it TSP-perfect}: in such graphs the TSP subtour formulation defines the convex hull of all vectors that satisfy the TSP requirements, and thus the solution to its linear relaxation is integer. The authors also showed that series parallel graphs are not TSP-perfect despite the GTSP being efficiently solvable in such graphs. \cite{fonlupt1993} extended these results by defining a generalisation of the TSP which, together with dynamic programming techniques, is used to show that the GTSP is polynomially solvable on a broader class of graphs (such as Halin graphs \citep{cornuejols1983}).

Unlike Ratliff \& Rosenthal, we assume that warehouses may be composed of multiple blocks. We however adopt the same representation as the corresponding graph is composed of series parallel subgraphs joined by paths of edges that represent cross-aisles. Figures \ref{fig2} and \ref{fig3} illustrate this as each block individually forms a single series parallel graph, however these two subgraphs are joined by the middle cross-aisle. Unlike Ratliff \& Rosenthal, we also assume that from the source one can reach any of the artificial vertices in the top cross-aisle. To the best of our knowledge, the GTSP in graphs composed of several connected series parallel subgraphs, such as warehouses with multiple blocks, has not been shown to be efficiently solvable.

A valid tour for a single picker includes vertices that must be visited and vertices that may be visited. This is due to the existence of artificial vertices and the fact that pickers in general are responsible for collecting a subset of orders which may not span the entire set of non-artificial vertices. For a single picker, the problem of finding the shortest walk which visits a known required subset of vertices is known as the Steiner TSP (STSP) and was first studied by \cite{cornuejols1985,fleischmann1985,orloff1974}. Several formulations for the STSP, including compact ones, were proposed by \cite{letchford2013}. The compact formulations introduced in \cite{valle2016} were inspired by similar formulations introduced in \citep{letchford2013}.

When trolley capacity and multiple pickers are considered, the problem is similar to the Steiner Capacitated VRP (SCVRP) \citep{letchford2006}. To the best of our knowledge no work in the literature has attempted to solve it directly, although many algorithms have been developed for similar problems such as the Capacitated VRP \citep{naddef2007, fukasawa2006} and the Capacitated Arc Routing Problem \citep{belenguer1998, letchford2009}. A closely related problem is the Mixed-Capacitated General Routing Problem (MCGRP), which reduces to the SCVRP if the set of required arcs is empty. \cite{bosco2013} and \cite{irnich2015} proposed branch-and-cut algorithms for the MCGRP, although the authors did not experiment on SCVRP instances (all instances had required arcs).

An important difference between SCVRP and JOBPRP however is that the indivisibility of orders must be respected; an order being assigned to a trolley implies that a group of vertices must be visited (assuming each item can be picked at just one vertex). SCVRP is a special case of JOBPRP whenever each order contains just a single product.

Many heuristics have been proposed for the problem of routing a single picker in warehouses with multiple blocks. An evaluation of several heuristics is given by \cite{petersen1997}. \cite{roodbergen2001_2} introduced a dynamic programming algorithm for picker routing in warehouses with up to three cross-aisles and \cite{roodbergen2001} presented heuristics for warehouses with multiple cross-aisles. The latter is compared to a branch-and-bound algorithm based on a classical TSP formulation of the problem, which assumes a  Hamiltonian circuit in a complete graph where all vertices must be visited exactly once (no artificial vertices are considered). \cite{theys2010} adapted the TSP LKH heuristic by \cite{lin1973} for single-picker routing and \cite{dekker2004} studied the case of a Dutch retail company whose layout is composed of multiple cross-aisles, two floors and different origin and destination points for the pickers. 

Several heuristics have also been proposed for batching and routing multiple pickers. Very often, routing distances are estimated using single-picker heuristics during the solution of the batching problem. \cite{henn2012_3} includes an extensive survey of batching methods.

For warehouses with two cross-aisles, a VNS heuristic was proposed by \cite{albareda2009} and Tabu Search and Hill Climbing heuristics were introduced by \cite{henn2012}. The latter methods were adapted by \cite{henn2013} to a problem that also considers order sequencing. \cite{azadnia2013} proposed heuristics that solve sequencing and batching first, then routing in a second stage. For warehouses with three cross-aisles, \cite{matusiak2014} proposed a simulated annealing algorithm which includes precedence constraints (for instance, a heavier box has to be at the bottom of a container); routing of candidate batches is solved with the A*-algorithm from \cite{hart1968}.

For general warehouses with any number of blocks, several batching methods are presented by \cite{deKoster1999}, where routing distances are estimated using single-picker methods. \cite{gibson1992} introduced batching heuristics based on two and four-dimensional space-filling curves. A genetic algorithm that jointly considers both routing and batching was proposed by \cite{hsu2004}. \cite{lam2014} introduced an integer programming formulation for batching, where routing distances are estimated. The problem is solved with a heuristic based on fuzzy logic. An association-based clustering approach to batching considering customer demands (instead of routing distance) is presented by \cite{chen2004}.

As far as we are aware, very few exact algorithms for the joint batching and routing problem have been proposed in the literature to date.  The work presented here builds upon our previous work \citep{valle2016}, as discussed elsewhere in this paper.  
\cite{hong2012a} presented a formulation that depends on explicitly enumerating all possible picker routes. No computational results for that formulation were presented, rather lower bounds on the solution to that formulation were developed as well as a heuristic procedure.
\cite{hong2012} dealt with multiple pickers and a single block (so just two cross-aisles). They considered aisle congestion where pickers are delayed by the presence of other pickers in the same aisle. In terms of picker routing they considered each picker to follow the same route around all of the aisles, with the proviso that some aisles can be skipped if desirable. 
\cite{hong2016} built on \cite{hong2012} and considered the situation when the pickers operate in a bucket brigade manner. In such a system each trolley (equivalently bucket) follows the same route around all of the aisles, being passed from picker to picker along the route.
\cite{scholz2016} dealt with just a single picker and a single block (so just two cross-aisles). Although they indicated that their approach can be extended to more than one block, no numerical results were given. Unlike these works, the exact algorithm proposed here is general in the sense that it can be applied to warehouses with any number of blocks.

\section{The joint problem of order batching and picker routing: a graph optimisation problem}
\label{sec:graph}

Let $\mathcal{P}$ be the set of products whose storage slots in the warehouse are known and $\mathcal{L}$ be the set of locations in the warehouse from where a picker can collect products, i.e. the middle of an aisle containing products on both sides, in different shelves, and from where the picker can reach those products. Each location $L \in \mathcal{L}$ contains a subset of products. A warehouse graph representation includes a vertex for every location, represented by the white vertices in Figure \ref{fig2}.

Let $O$ be the set of orders that must be collected. Each order $o \in O$ contains a subset of products $P_o \subseteq \mathcal{P}$. Accordingly, for each order $o \in O$, the subset of locations $L_o \subseteq \mathcal{L}$ contains all products in $P_o$ - it is possible that multiple products of the same order are in a single location. 
$L(O) = \bigcup_{o \in O} L_o$ represents the set of locations that contains all products that need to be picked in all orders. Also let $d_{lm} \geq 0$ be the distance between locations $l, m \in \mathcal{L}$, symmetric so that $d_{lm} = d_{ml}$.

In practice, very often $L(O) \not = \mathcal{L}$, i.e. there are locations for which none of their products are present in any of the outstanding orders. In such cases, the graph that represents a warehouse can be reduced by eliminating all vertices that represent locations in $\mathcal{L} \setminus L(O)$. An example of such reduced graph can be seen in Figure \ref{fig3}, where several vertices were removed as discussed previously. If $v$ is eliminated, and there are arcs $(u,v), (v,u), (w,v), (v,w)$, we create arcs $(u,w), (w,u)$ where $d_{uw} = d_{wu} = d_{uv} + d_{vw} = d_{wv} + d_{vu}$.

Let ${\cal T}$ denote the set of available pickers (or trolleys), $T = |{\cal T}|$, $B$ be the number of baskets that a trolley can carry and $b_o$ be the known number of baskets needed to carry order $o \in O$. We assume that a basket will only contain products from a single order, even if it is partially empty, i.e., it is not possible to put products of different orders in the same basket. Finally, let $s$ be the origin point from where the trolleys depart (and to where they must return), and accordingly let $d_{sl} \geq 0$ be the distance between location $l$ and $s$.

To define JOBPRP, we introduce a directed and connected graph $D = (V,A)$. The set of vertices $V$ is given by the union of $s$ (the origin point), a set $V(O)$ containing a vertex for every location $l \in L(O)$ and a set $V_A$ of ``artificial locations''. These artificial locations are placed in corners between aisles and cross-aisles and do not contain products that must be picked. We also define sets $V_o$ containing a vertex for every location $l \in L_o$, $\bigcup_{o \in O} V_o = V(O)$. Thus, $V = \{s\} \cup V(O) \cup V_A$ and $|V| = 1 + |V(O)| + |V_A|$. In Figure \ref{fig3}, the origin represents $s$, the 17 white vertices represent $V(O)$ and the 9 black vertices represent $V_A$. 

Given a vertex $i \in V(O)$, define $l(i): V(O) \rightarrow L(O)$ as the location to which the vertex refers. Arcs in set $A$ connect two neighbour vertices of $V$. For instance, arc $(i,j)$ connects location $l(i)$ to location $l(j)$. $d_{l(i),l(j)} \geq 0$ represents the distance between vertices $i$ and $j$. For readability we write $d_ {ij}$ instead of $d_{l(i),l(j)}$.

A solution to JOBPRP in $D$ is a collection of $T^* \leq T$ closed walks, one for each trolley. Each closed walk $t = 1,\dots, T^*$ is associated with a subgraph $(V_t \cup \{s\},A_t)$ of $D$.
Each walk starts at $s$, traverses a set $A_t \subseteq A$ of selected arcs and returns to $s$. 
For each $t = 1,\dots, T^*$, let $O_t \subseteq O$ denote the subset of orders collected during walk $t$.


The following requirements must be met by solutions to JOBPRP: 
\begin{itemize}
\item  Capacity constraints impose that, for every $t = 1,\dots,T^*$,  $\sum_{o \in O_t} b_o \leq B$. 

\item Indivisibility of the products $P_o$ in each order $o \in O_t$ collected by each trolley imposes that, for every $o \in O_t$, there exists at least one walk $t$ where if $i \in V_t: l(i) \in L_o$, then every $j: l(j) \in L_o$ must also be visited by $t$, i.e., $j \in V_t$. 

\item Finally, in order to guarantee that the $T^*$ walks collect all orders, we must impose that $\cup_{t = 1}^{T^*} O_t = O$.  We assume that there are at least as many orders as there are trolleys available.
\end{itemize} 

JOBPRP is then the problem of finding $T^* \leq T$ closed walks, meeting the requirements outlined above, in order to minimize the total length $\sum_{t = 1}^{T^*} \sum_{(i,j) \in A_t} d_{ij}$.

As $D$ is modelled as a sparse graph, vertices in $V(O)$ have only two neighbours while those in $V_A$ do not connect to more than four other vertices. A feasible trolley walk for JOBPRP may be forced to not only repeat vertices but also to visit vertices from non-collected orders. For instance, this could occur if a product from a non-collected order happens to be between two distinct products in a collected order. However, under reasonable assumptions, we can prove that the optimal solution to JOBPRP is composed of walks that may repeat vertices, but will not repeat arcs. An explicit proof can be found in the appendix of \cite{letchford2009}, for a closely related problem. 

\begin{lemma}
Let us assume that $D$ is such that for every arc $(i,j)$ there is a corresponding arc $(j,i)$ and that $d_{ij} = d_{ji} \geq 0$. Then, in the optimal solution $T^*$ to JOBPRP, every walk $t \in T^*$ visits each arc in $A_t$ exactly once.
\label{lemma1}
\end{lemma}

\section{An integer programming formulation for the JOBPRP}
\label{sec:IP}

Based on Lemma \ref{lemma1}, we introduce an integer programming formulation for the JOBPRP. The formulation allows vertices to be visited multiple times, but enforces the condition that each arc cannot be traversed more than once. The model uses exponentially many constraints to enforce connectivity of the closed walks. 

The formulation uses binary decision variables $z_{ot}$ to indicate whether ($z_{ot}=1$) or not ($z_{ot}=0$) trolley $t$ picks order $o \in O$, $x_{tij}$ to indicate whether ($x_{tij} = 1$) or not ($x_{tij} = 0$) arc $(i,j) \in A$ is traversed by trolley $t$, $\alpha_{t}$ to indicate whether ($\alpha_{t} = 1$) or not ($\alpha_{t} = 0$) trolley $t$ picks at least one order and $y_{ti}$ to indicate whether ($y_{ti} = 1$) or not ($y_{ti} = 0$) vertex $i \in V \setminus \{s\}$ is visited by trolley $t$. In addition, the model uses variables $g_{ti} \in \mathbb{Z}_+$ to indicate the outdegree of vertex $i \in V$ in the closed walk for trolley $t$.

Let $\delta^-(W) = \{ (i,j) \in A: i \not \in W, j \in W\}$, $\delta^+(W) = \{ (i,j) \in A: i  \in W, j \not \in W\}$ and $A(W) = \{ (i,j) \in A: i  \in W, j \in W\}$. Also, let $|\delta^+(i)|$ be the maximum outdegree of $i \in V$ (for simplicity we write $\delta^+(i)$ instead of $\delta^+(\{i\})$). JOBPRP can be stated as the following Integer Program:

\begin{equation}
\min  \sum_{t \in \cal T} \sum_{(i,j) \in A} d_{ij} x_{tij}
\label{obj1}
\end{equation}

\noindent subject to:

\begin{optprog}
& \sum_{o \in O} b_o z_{ot} & \;\leq \;& B \alpha_{t}, & \forall t \in {\cal T} \label{eq1}\\

& \sum_{t \in \cal T} z_{ot} & \;= \;& 1, & \forall o \in O \label{eq2}\\

& \sum_{(i,j) \in \delta^+(i)} \kern -1em x_{tij} & \;\geq\; & z_{ot}, & \forall o \in O, t \in {\cal T}, i : l(i) \in L_o,  \label{eq3}\\

& \sum_{(i,j) \in \delta^+(i)} \kern -1em x_{tij} & \;= \;& \kern -1em \sum_{(j,i) \in \delta^-(i)} \kern -1em x_{tji}, & \forall i \in V, t \in {\cal T} \label{eq4}\\

& \sum_{(s,j) \in \delta^+(s)} \kern -1em x_{tsj} = \kern -1em \sum_{(j,s) \in \delta^-(s)} \kern -1em x_{tjs} & \; = \;& \alpha_{t}, &  \forall t \in {\cal T} \label{eq5} \\

& x_{tij} & \leq & \alpha_t, & \forall (i,j) \in A, t \in {\cal T} \label{eq6} \\

& z_{ot} & \leq & \alpha_t, & \forall o \in O, t \in {\cal T} \label{eq7} \\

& \sum_{o \in O} z_{ot} & \geq & \alpha_t, & t \in {\cal T} \label{eq8} \\

& \sum_{(i,j) \in \delta^+(i)} \kern -1em x_{tij} & \;= \;& g_{ti}, & \forall i \in V, t \in {\cal T} \label{eq9}\\

& y_{ti}                               & \;\geq \;& x_{tij}, & \forall (i,j) \in A, t \in {\cal T} \label{eq10} \\

& \sum_{j \in W} g_{tj} & \;\geq \;& \;y_{ti} + \kern -1em \sum_{(j,k) \in A(W)} \kern -1em x_{tjk} , & \forall i \in W, W \subseteq V\setminus \{s\}, |W| > 1, t \in {\cal T} \label{eq11} \\

& x_{tij} & \in & \mathbb{B}, & \forall (i,j) \in A, t \in {\cal T} \label{eq12}\\

& z_{ot}  & \in & \mathbb{B}, & \forall o \in O, t \in {\cal T} \label{eq13}\\

& 0 \leq \alpha_{t}  & \leq & 1, & t \in {\cal T} \label{eq14} \\
%

& y_{ti}  & \leq & \alpha_t, & \forall i \in V, t \in {\cal T} \label{eq15}
\end{optprog}

Constraints (\ref{eq1}) guarantee that the number of baskets in a trolley does not exceed its available capacity while constraints (\ref{eq2}) ensure that each order is collected by precisely one trolley. 
Constraints (\ref{eq3}) enforce that if an order is assigned to a trolley, then some vertex that stores a product of this order will be visited by the trolley at least once. 
Constraints (\ref{eq4}) make sure that for every arc that reaches a vertex, there is one that leaves it. Constraints (\ref{eq5}) ensure that if a trolley picks any order, then it necessarily departs from the origin. Without significant loss of generality we assume that the trolley visits the source only once.
Constraints (\ref{eq6}) and (\ref{eq7}) make sure that a trolley visits an arc or picks an order only if it is used, while constraints (\ref{eq8}) ensure that if a trolley is used, then at least one order is picked by that trolley. 
Constraints (\ref{eq9}) and (\ref{eq10}) define the outdegree and the $y$ variables for each vertex. Note that if the maximum outdegree of each vertex was not allowed to be greater than one, but instead, if $g_{ti} \in \{0,1\}$, we would have $y_{ti} = g_{ti}$, and constraints (\ref{eq11}) would change to the generalized subtour breaking constraints $\sum_{(j,k) \in A(W)}x_{tjk} \leq \sum_{j \in W \setminus \{i\}}y_{tj}$ \citep{lucena1992}. 
Constraints (\ref{eq11}) do allow subtours found in closed walks as long as at least one vertex in the cycle has an outdegree of $2$.
Constraints (\ref{eq12})-(\ref{eq15}) deal with the variables. Note here that although in constraint (\ref{eq14}) we allow the $\alpha_t$ variables to be fractional they will be forced to be binary due to the influence of other constraints (e.g.~(\ref{eq12}) and (\ref{eq5})). 



\subsection{Symmetry breaking constraints}

Due to the artificial indexation of identically sized trolleys, formulation (\ref{eq1})-(\ref{eq15}) above suffers significantly from symmetry. That means that identical order to walk assignments lead to different vector representations, when we simply change the indices of the trolleys to which the closed walks are assigned. Branch-and-bound algorithms based on symmetric formulations tend to perform poorly, since they enumerate search regions that essentially lead to the same solution. As an attempt to overcome symmetry for JOBPRP, we add the following constraints:

\begin{equation}
\sum_{t = 1}^{o} z_{ot} \geq 1 ,\;\;\;\;\; o = 1, \ldots, T
\label{eq16}
\end{equation}

\noindent Here we ensure that the first order, $o = 1$, will be picked by the trolley assigned to the first index $t = 1$. The second order, $o = 2$, is either picked by trolley 1 or 2, and so forth. 

We also add the following constraints:

\begin{equation}
\alpha_{t} = 1, \;\;\;\;\;   t = 1, \ldots, \Big\lceil \sum_{o \in O} b_o / B \Big\rceil 
\label{eq17}
\end{equation}

\noindent The expression $\Big\lceil \sum_{o \in O} b_o / B \Big\rceil$ represents the known minimum number of trolleys necessary, constraints (\ref{eq17}) enforce that the first trolleys in ascending index order are selected.

\subsection{Implementation details}
\label{sec:impl}

In this section, we highlight some implementation details regarding our solution of the formulation outlined above.

\subsubsection{Branch-and-cut}
\label{sec:branch}

We employ a branch-and-cut algorithm \citep{padberg1987, padberg1991} which separates inequalities (\ref{eq11}). For candidate integral solutions, we verify via a depth-first search whether the graph is connected, if not we add one constraint of type (\ref{eq11}) for every disconnected component. 

Optionally, we may also check the existence of violated constraints (\ref{eq11}) in connected fractional solutions. Let $\overline{g}_{ti}$, $\overline{y}_{ti}$ and $\overline{x}_{tij}$ be the values taken by the corresponding variables $g_{ti}$, $y_{ti}$ and $x_{tij}$ in a linear programming relaxation of formulation (\ref{eq1})-(\ref{eq17}). For every $t = 1, \ldots, T$, let $\overline{V}^t \subseteq V \setminus s$ be the set of vertices where $\overline{y}_{ti} > 0$ and $A(\overline{V}^t) \subseteq A$ be the set of arcs with both ends in $\overline{V}^t$ where $\overline{x}_{tij} > 0$. The problem consists of finding a subset of vertices $\overline{W} \subseteq \overline{V}^t$ and a vertex $i \in \overline{W}$ for which a constraint of type (\ref{eq11}) is violated. 

This problem can be polynomially solved by finding the minimum cut (maximum flow) on a network given by the graph $\overline{D}^t = (\overline{V}^t$, $A(\overline{V}^t))$ with capacities $\overline{x}_{tij} \forall (i,j) \in A(\overline{V}^t)$. Given any arbitrary $\overline{W} \subset \overline{V}^t$, constraints (\ref{eq10}) ensure that $\sum_{j \in {\overline{W}}} \overline{g}_{tj} - \sum_{(j,k) \in A(\overline{W})} \overline{x}_{tjk} = \sum_{(j,l) \in \delta^+(\overline{W})}  \overline{x}_{tjl}$, which is exactly the value of the cut that separates $\overline{W}$ from $\overline{V}^t \setminus \overline{W}$.

To solve the minimum cut/maximum flow problem, we employ the push-relabel algorithm introduced by \cite{goldberg1988}. Since the problem has to be solved separately for every $t = 1,\ldots,T$, a violated cut for $t$ is temporarily kept in a cut pool and checked for violation in subsequent iterations for other trolleys.

Whenever the algorithm is configured to only check candidate integral solutions for violations, we use the term Integral Branch-and-Cut, or {\it \textbf{IBC}} for short. Whenever violations are checked for both integral and fractional solutions, we use the term Fractional Branch-and-Cut, or {\it \textbf{FBC}}.

\subsubsection{Heuristic}
\label{sec:heuristic}

We implemented a constructive heuristic to provide an initial solution. Batching is computed via the time savings heuristic as explained in Section 2.2.1 of \cite{deKoster1999}, but the savings matrix is not recalculated as orders are clustered. 

For the estimation of partial route distances, we employed picker routing heuristics introduced in \cite{roodbergen2001}, which are specific for warehouses with multiple cross-aisles. Partial routes are computed by running each of the following heuristics: s-shape, largest gap, combined and combined+. The best solution among these four different methods is taken as the input value in the savings matrix.

\subsubsection{Solver tuning} 

We employed CPLEX 12.6.0 \citep{cplex} as the mixed-integer programming solver. Branching priority is given to $z_{ot}$ variables. CPLEX presolve is enabled but limited due to callback calls within the optimisation. 

CPLEX has built-in general-purpose separation algorithms for certain classes of inequalities. After some computational experimentation we observed that some of these inequalities help to improve linear relaxation bounds; on the other hand more computation time is required to solve each node in the branch-and-bound tree (with a greater emphasis on the first node). We disabled most CPLEX separation procedures due to them being ineffective, however we set the search for violated lift-and-project, zero-half, mixed-integer rounding, Gomory and cover cuts to aggressive mode. This improved results slightly, however smaller instances can actually be solved faster by disabling all separation algorithms.

All other solver settings are left as default.

\section{Valid inequalities}
\label{sec:ineq}

\begin{figure}[!htb]
\centering
  \includegraphics[width=0.5\textwidth]{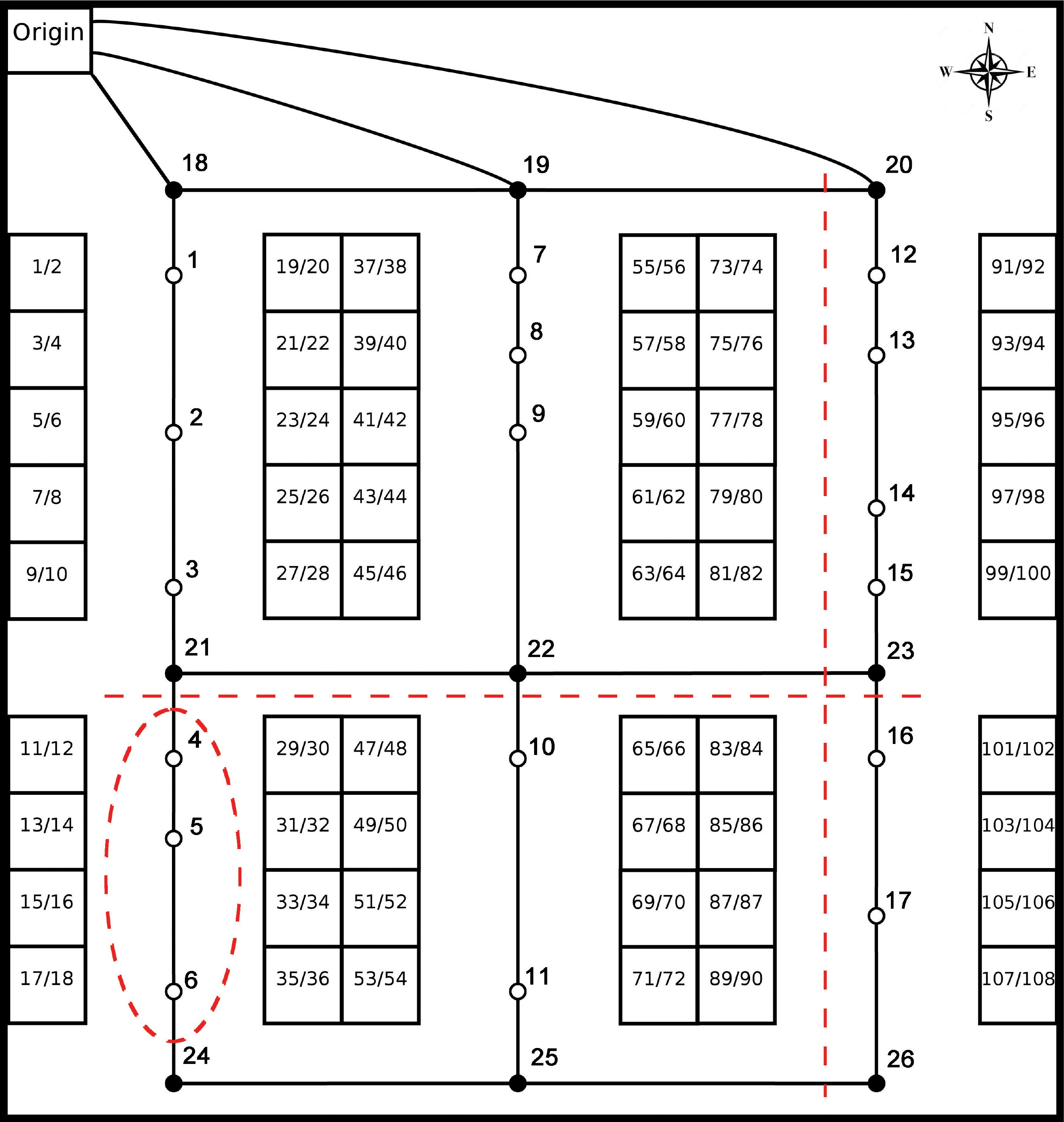}
  \caption{Warehouse example}
  \label{fig4}
\end{figure}

In this section we introduce classes of JOBPRP valid inequalities based on the standard layout of warehouses. All these inequalities are valid for formulation (\ref{eq1})-(\ref{eq17}) as well as for the two compact formulations introduced in \cite{valle2016}. 

With regard to a terminological issue here we refer below to \emph{\enquote{valid inequalities}}. We have used this phrase in its most general form, by which we mean that we have presented inequalities which are valid when we are seeking the optimal solution to the integer program that we have formulated. In polyhedral studies the term \enquote{valid inequalities} often has a more restricted meaning and refers to linear inequalities that are satisfied by all feasible integer solutions. Other authors also use the term \emph{\enquote{optimality cuts}} or \emph{\enquote{logic cuts}} \citep{hooker1994} for inequalities that cut off sub-optimal integer feasible solutions.

Let $\mathbb{A}$ be the set of aisles with $W_A=|\mathbb{A}|$ and let
$\mathbb{C}$ be the set of cross-aisles with  $W_C=|\mathbb{C}|$. 
 Aisles are numbered $[1, \ldots, W_A]$ from left to right, and cross-aisles are numbered $[1, \ldots, W_C]$ from top to bottom. The number of subaisles is equal to $W_A(W_C-1)$.  We assume without significant loss of generality that the source, from where trolleys depart and to which they return, is located at the top left corner of the warehouse.

For ease of illustration of the valid inequalities introduced below we shall make use of Figure~\ref{fig4}. That figure contains three aisles and three cross-aisles, with six subaisles. Here the vertices correspond to the reduced graph (cf~Figure~\ref{fig3} seen before) corresponding to the orders that must be picked. The dotted lines seen in Figure~\ref{fig4} correspond to cuts considered below. For ease of discussion we shall use compass directions (north, south, east, west) as also shown in Figure~\ref{fig4}. 

With respect to the labelling of the vertices seen in Figure~\ref{fig4} we have labelled the vertices in each cross-aisle and each subaisle  in consecutive order. Let the artificial vertex at the intersection of aisle $a$ and cross-aisle $c$ be $v(a,c)$. Let the ordered set of vertices in the subaisle associated with aisle $a$ that lies between cross-aisle $c$ and cross-aisle $c+1$ be $S(a,c)$. Referring to Figure~\ref{fig4} for example and taking the subaisle with a dotted ellipse in it (so $a=1$ and $c=2$) we have $v(1,2)=21$, $v(1,3)=24$ and $S(1,2)=\{4,5,6\}$. Let $R(a,c)$ denote the cardinality of set $S(a,c)$, so $R(a,c)$ is the number of vertices in the subaisle below the intersection of aisle $a$ and cross-aisle $c$. We use the notation $S(a,c,r)$ to refer to the $r$'th element of the ordered set $S(a,c)$ (so $S(a,c,r)$ is the $r$'th vertex in the subaisle below  the intersection of aisle $a$ and cross-aisle $c$). To ease the presentation below we define $S(a,c,0)=v(a,c)$ and $S(a,c,R(a,c)+1)=v(a,c+1)$. In this section, whenever convenient, we express a variable $x_{tij}$ as $x_{t,i,j}$.

We now consider the valid inequalities (cuts) that were found to lead to improved computational performance. 

\subsection{Further symmetry breaking}

The additional notation allows us to introduce a further symmetry-breaking constraint. In order to break symmetry in the directions adopted by each trolley walk then, as distance is a symmetric function, we can insist that the arc out of $s$ for a trolley is to the left (so west) of the arc into $s$. For example, considering Figure~\ref{fig4} if we use the arc $(s,19)$ for a trolley then we can only use arc $(19,s)$ or $(20,s)$ for the return to $s$. The constraints which enforce this condition are:

\begin{equation}
\sum_{k=a}^{W_A} x_{t,v(k,1),s} \geq \sum_{k=a}^{W_A} x_{t,s,v(k,1)}, \;\;\;\;\;\; \forall a \in \mathbb{A}\setminus[1], t \in {\cal T} 
\label{jeb1}
\end{equation}

Although technically the constraint above could be intuitively represented as $\sum_{k=a}^{W_A} x_{t,v(k,1),s} \geq x_{t,s,v(a,1)}$, constraints in the form (\ref{jeb1}) are stronger since decision variables are nonnegative. 

\subsection{First cross-aisle vertex}

Consider a trolley that leaves the origin $s$ in Figure~\ref{fig4}.   It must go to one of the artificial  vertices in the first cross-aisle. Now if it goes to one of these vertices it makes no sense for that trolley to next go to an adjacent artificial  vertex in the same cross-aisle since (under reasonable assumptions about the distance function) it would be of lesser distance to go directly from $s$ to that adjacent vertex. We can therefore deduce that from the first artificial  vertex visited in the first cross-aisle the trolley proceeds down the associated subaisle in a north to south direction. The constraints which enforce this are:

\begin{equation}
x_{t,v(a,1),S(a,1,1)} \geq x_{t,s,v(a,1)}, \;\;\;\;\;\; \forall a \in \mathbb{A}, t \in {\cal T} 
\label{jeb2}
\end{equation}

Using a similar argument it is trivial to show that at the end of a trolley walk, when the trolley returns to the origin $s$   from an artificial vertex in the first cross-aisle, it must have come to that artificial vertex from the vertex below it in the associated subaisle. The constraints which enforce this are:
\begin{equation}
x_{t,S(a,1,1),v(a,1) } \geq x_{t,v(a,1),s}, \;\;\;\;\;\; \forall a \in \mathbb{A}, t \in {\cal T} 
\label{jeb3}
\end{equation}

In fact we can go further here with respect to the first subaisle entered after leaving the origin. Suppose in Figure~\ref{fig4} a trolley leaves the origin, goes to vertex 19 and enters the subaisle below vertex 19. Potentially it could visit all of the vertices in that subaisle (vertices 7, 8, 9) and then return to the origin using arc $(19,s)$.   However if there are any items in orders assigned to trolley $t$ that cannot be completely picked from the subaisle then the trolley must either use (9,22) to reach the cross-aisle and visit other aisles, or it can reverse back to 19 but must then use one of (19,20) or (19,18) to visit other aisles. In other words $x_{t,9,22} + x_{t,19,18} + x_{t,19,20}  \geq x_{t,s,19} -(1-z_{ot})$ for all orders $o$ that cannot be completely satisfied using the subaisle alone.

Generalising, let $\Phi(a)$ be the set of orders $o \in O$ that cannot be completely satisfied using all of the vertices in the subaisle between artificial vertex $v(a,1)$  at the intersection of the first cross-aisle  and aisle $a$ and $v(a,2)$, so the vertices in $S(a,1)$. Then we have:
\begin{equation}
\begin{split}
x_{t,S(a,1,R(a,1)),v(a,2)}  +  x_{t,v(a,1),v(a-1,1)} + x_{t,v(a,1),v(a+1,1)} \geq x_{t,s,v(a,1)} - (1-z_{ot}), \\
\;\;\;\;\;\; \forall o \in \Phi(a), a \in \mathbb{A}, t \in {\cal T}
\label{jeb4}
\end{split}
\end{equation}

Here, for notational convenience, we disregard any terms (such as occur when $a=1$ or $a=W_A$) that have no meaning with respect to the variables that we have defined.

\subsection{Cross-aisle cuts}

Consider Figure~\ref{fig4}  where we show a cross-aisle cut as dotted below vertices 21, 22 and 23. That cross-aisle cut separates the supermarket into two portions, one portion containing the origin north of the cut, the other portion containing some subaisles south of the cut. 

Now suppose for some order $o$ that there are items in the order which only lie in the lower portion (more technically can only be picked using vertices in the lower southern portion). Then if trolley $t$ picks order $o$ it must cross the cut at least once in a north to south direction. For the cut shown the constraint which enforces this would be $x_{t,21,4}+x_{t,22,10}+x_{t,23,16} \geq z_{ot}$. Obviously the trolley must also cross the cut in a south to north direction an equal number of times 
(i.e.~we also have $x_{t,4,21}+x_{t,10,22}+x_{t,16,23}=
x_{t,21,4}+x_{t,22,10}+x_{t,23,16}$).  

Define the set $\Gamma(c)$ as the set of orders $o \in O$ for which there is at least one item $i$  that has a unique picking location at a vertex in the southern portion when a cross-aisle cut is made just below cross-aisle $c$ ($c=1,\ldots,W_C-1$). It is clear from Figure~\ref{fig4}  that  $\Gamma(1)$ will contain all orders. 
Then the cross-aisle cut constraints which enforce the condition that the cut is crossed from north to south at least once, and from south to north an equal number of times, are:
\begin{optprog}
& \sum_{a \in \mathbb{A}}x_{t,v(a,c),S(a,c,1)} & \;\geq\; & z_{ot}, 
&\;\;\;\;\;\;  \forall c \in \mathbb{C}\setminus[W_C], o \in \Gamma(c), t \in {\cal T} \label{jeb5} \\

& \sum_{a \in \mathbb{A}} x_{t,S(a,c,1),v(a,c)} & \;=\; & \sum_{a \in \mathbb{A}}x_{t,v(a,c),S(a,c,1)},  
&\;\;\;\;\;\;\; \forall c \in \mathbb{C}\setminus[W_C],   t \in {\cal T}  \label{jeb6}
\end{optprog}

As well as these constraints we also have constraints associated with positioning the cross-aisle cut just above the artificial vertices in a cross-aisle. For example, if we have a cross-aisle cut positioned just above the second cross-aisle in Figure~\ref{fig4}, then the constraint that ensures that this cut is crossed would be $x_{t,3,21}+x_{t,9,22}+x_{t,15,23} \geq z_{ot}$.

The reason why a constraint above the cross-aisle  is not redundant when we also have a constraint below the cross-aisle can be illustrated using Figure~\ref{fig4}.
In the absence of the complete set of subtour elimination constraints it is possible to have a trolley subtour around the four artificial vertices 21, 22, 25, 24 that would satisfy the cross-aisle cut positioned just below the second cross-aisle, but would not satisfy a cross-aisle cut positioned just above the second cross-aisle.

Then the cross-aisle cut constraints which enforce the condition that a cut positioned just above a cross-aisle 
is crossed from north to south at least once, and from south to north an equal number of times, are:

\begin{optprog}
& \sum_{a \in \mathbb{A}}x_{t,S(a,c-1,R(a,c-1)),v(a,c)} & \;\geq\; & z_{ot}, 
&\;\; \forall c \in \mathbb{C}\setminus[1], o \in \Gamma(c), t \in {\cal T} \label{jeb7} \\

& \sum_{a \in \mathbb{A}}x_{t,v(a,c),S(a,c-1,R(a,c-1))} & \;=\; & \sum_{a \in \mathbb{A}}x_{t,S(a,c-1,R(a,c-1)),v(a,c)}, 
&\;\; \forall c \in \mathbb{C}\setminus[1],  t \in {\cal T}  
 \label{jeb8}
\end{optprog}

\subsection{Aisle cuts}

Consider Figure~\ref{fig4} and the aisle cut shown there dotted to the left of the vertices in the third aisle. That aisle cut separates the supermarket into two portions, one western portion containing the first two aisles and an eastern portion containing the third aisle. 

Suppose that for some order $o$ there is at least one item $i$ (in order $o$) with a unique location for picking in the first two aisles, and in addition there is at least one item $j$ (in order $o$) with a unique location for picking in the third aisle. Then if trolley $t$ picks order $o$ the trolley must cross the cut (either from west to east or from east to west). For the aisle cut shown the constraint that enforces this is $x_{t,19,20}+ x_{t,20,19}+x_{t,22,23}+ x_{t,23,22}+x_{t,25,26} + x_{t,26,25} \geq z_{ot}$.

This aisle cut has a similar logic to the cross-aisle cut considered above, except that here the situation is more complicated as we have the possibility of entering the aisle directly from the origin (and indeed returning directly to the origin from the same aisle). Note here that
for this cut to be valid it has to be true that the trolley is not allowed to revisit the origin  (i.e.~we exclude from the walk the case: trolley to item $i$ in the first two aisles, then back to $s$, then from $s$  to the artificial  vertex at the end of the third aisle, then to item $j$ in the third aisle and then back to $s$). Equation (\ref{eq5}) above ensures  this.

Generalising, consider an aisle cut in aisle $a$. Define the set $\Theta(a)$ as that set of orders $o \in O$ such there is at least one item $i$ (in order $o$) with a unique location for picking in aisles $1,\ldots,a-1$; and in addition there is at least one item $j$ (in order $o$) with a unique location for picking in aisles $a,\ldots,W_A$. Then the aisle cut constraint is:
\begin{equation}
\sum_{c \in \mathbb{C}} (x_{t,v(a-1,c),v(a,c)} + x_{t,v(a,c),v(a-1,c)}) \geq  z_{ot},
\;\;\;\;\;\; \forall a \in \mathbb{A}\setminus[1], o \in \Theta(a), t \in {\cal T} 
\label{jeb9}
\end{equation}

\subsection{Subaisle cuts}
Consider the subaisle in Figure~\ref{fig4} where we have vertices 4,5 and 6 within a dotted ellipse. Now suppose that vertex 5 within this subaisle is the unique picking location for some item contained in order $o$. Then if trolley $t$ picks order $o$  it must be the case  that trolley $t$ enters the subaisle, either from the north so using the arcs (21,4) and (4,5) to reach vertex 5, or from the south  so using the arcs (24,6) and (6,5) to reach vertex 5. 


Let $M^N(t,q)$ be the minimum $x_{tij}$ value over all arcs $(i,j)$ in the unique path from the northern cross-aisle artificial vertex to a vertex $q$ in the subaisle associated with the cross-aisle artificial vertex. Let $M^S(t,q)$ be the minimum $x_{tij}$ value over all arcs in the unique path from the southern cross-aisle artificial vertex to vertex $q$. 

So informally for vertex 5 in Figure~\ref{fig4}  these values are the minimum arc $x$ values when coming from the north to vertex 5 (so coming from artificial vertex 21) or coming from the south to vertex 5 (so coming from artificial vertex 24). For vertex 5 the northern $M$ values can be defined using:
$ M^N(t,4) = x_{t,21,4}$;   
$ M^N(t,5)  \leq   x_{t,4,5}$ and 
$ M^N(t,5)  \leq  M^N(t,4)$.
Here the first equation initialises the $M$ value for the first vertex (vertex 4) in the path north to south in the subaisle, the second and third constraints ensure that $M$ value associated with vertex 5 is bounded appropriately.
The southern $M$ values for vertex 5 can be defined in a similar way using:
$ M^S(t,6) =  x_{t,24,6}$; 
$ M^S(t,5) \leq   x_{t,6,5}$ and 
$ M^S(t,5) \leq   M^S(t,6)$.

Then assuming (as above) that vertex 5 must be visited by trolley $t$ if order $o$ is picked by trolley $t$ we must have:
\begin{equation}
M^N(t,5) +  M^S(t,5) \geq z_{ot}  \label{jeb10}
\end{equation}

Clearly if we knew in advance whether the trolley to vertex 5 came from the north or from the south we could strengthen this constraint, equation~(\ref{jeb10}).  Although we do not know from which direction the trolley visits vertex 5 it must  however be true that, whatever direction the trolley comes from, the sum of the $M$ variables must be at least $z_{ot}$, which is one if trolley $t$ picks order $o$.

Although these constraints might appear less than intuitive at first sight they arose because we had clear examples where the solution to the LP relaxation of the problem violated these constraints. In relation to a minor technical issue here, note that although we have spoken above as to these $M$ values being the minimum values over an arc path, it is possible that the numeric values assigned by the optimiser to the $M$ variables could be less than the minimum value over a path. However, since constraint (\ref{jeb10}) is a $\geq$ constraint then alternative optimal solutions exist where the $M$ values are artificially increased to the minimum value over the path.

Generalising, consider the subaisle that lies between artificial vertices $v(a,c)$ and $v(a,c+1)$. This subaisle contains vertices $S(a,c,1),S(a,c,2),\ldots, S(a,c,R(a,c))$. Then the constraints that define the $M$ variables for this subaisle are as below, where these constraints apply for all such subaisles (so for
$\forall a \in \mathbb{A}$, $\forall c \in \mathbb{C}\setminus[W_C]$).
\begin{optprog}

\kern -0.5cm & M^N(t,S(a,c,1)) & \;=\; & x_{t,v(a,c),S(a,c,1)}  
& & \kern -0.2cm \forall a \in \mathbb{A}, c \in \mathbb{C}\setminus[W_C], t \in {\cal T}
\label{jeb11} \\

\kern -0.5cm & M^N(t,S(a,c,r)) & \;\leq\; &  x_{t,S(a,c,r-1),S(a,c,r)} 
& \kern -0.2cm r=2,\ldots,R(a,c), & \kern -0.2cm \forall a \in \mathbb{A}, c \in \mathbb{C}\setminus[W_C], t \in {\cal T}
\label{jeb12} \\

\kern -0.5cm & M^N(t,S(a,c,r)) & \;\leq\; &  M^N(t,S(a,c,r-1))  
& \kern -0.2cm r=2,\ldots,R(a,c), & \kern -0.2cm \forall a \in \mathbb{A}, c \in \mathbb{C}\setminus[W_C], t \in {\cal T}
\label{jeb13} \\

\kern -0.5cm & M^S(t,S(a,c,R(a,c))) & \;=\; & x_{t,v(a,c+1),S(a,c,R(a,c))} 
&  & \kern -0.2cm \forall a \in \mathbb{A}, c \in \mathbb{C}\setminus[W_C], t \in {\cal T}
 \label{jeb14} \\

\kern -0.5cm & M^S(t,S(a,c,r-1)) & \;\leq\; &  x_{t,S(a,c,r),S(a,c,r-1)} 
& \kern -0.2cm r=2,\ldots,R(a,c), & \kern -0.2cm \forall a \in \mathbb{A}, c \in \mathbb{C}\setminus[W_C], t \in {\cal T}
\label{jeb15} \\

\kern -0.5cm & M^S(t,S(a,c,r-1)) & \;\leq\; &  M^S(t,S(a,c,r))  
& \kern -0.2cm r=2,\ldots,R(a,c), & \kern -0.2cm \forall a \in \mathbb{A}, c \in \mathbb{C}\setminus[W_C], t \in {\cal T}
\label{jeb16}
\end{optprog}

Define the set $\Omega(p)$ as that set of orders $o \in O$ such there is at least one item in order $o$ with a unique location for picking at vertex $p \in V(O)$. Then the constraint relating the $M$ values for vertex $p$ to $z_{ot}$ is:

\begin{equation}
M^N(t,p) +  M^S(t,p) \geq z_{ot}, \;\;\;\;\;\; \forall p \in V(O), o \in \Omega(p), t \in {\cal T}  
\label{jeb17}
\end{equation}

Including these constraints involves an additional $2T|V(O)|$ variables and an additional $5T|V(O)|$ constraints (since essentially each $M$ variable is involved in two constraints in equations~(\ref{jeb11})-(\ref{jeb16})) as well as involved in equation~(\ref{jeb17}).

\subsection{Artificial vertex reversal}
\label{sec:artificialreversal}

Suppose  in Figure~\ref{fig4} we, for some trolley $t$,  use the arc (22,23) and once at vertex 23 simply reverse, so use the arc (23,22), but do not use either of the arcs (23,15)  and (23,16). The reason why this might occur is that reversing at artificial vertex 23 ensures that  the aisle cut associated with the third aisle (which here is $x_{t,19,20}+ x_{t,20,19}+x_{t,22,23}+ x_{t,23,22}+x_{t,25,26} + x_{t,26,25} \geq z_{ot}$) is satisfied. Clearly going from 22 to 23 and straight back again cannot be part of any optimal solution (since we are travelling extra distance for no practical effect).

In order to prevent this situation occurring we can impose additional constraints. For vertex 23  an appropriate constraint would be $x_{t,23,15}+ x_{t,23,16} \geq x_{t,22,23}$  which ensures that if we use arc (22,23) we must also go from vertex 23 to vertex 15 or vertex 16. Generalising, we have that: 

\begin{equation}
\sum_{\substack{j \in\delta^+(k) \\ j \neq i} } \kern -0.4em x_{tkj} \geq x_{tik} \;\;\;\;\;\; \forall k \in \delta^+(i),\;\; i,k \in \cup_{a \in \mathbb{A}} \cup_{c \in \mathbb{C}} \{v(a,c)\}, t \in {\cal T} 
\label{jeb18} 
\end{equation}

This constraint ensures that if for trolley $t$ we use arc $(i,k)$  between two artificial vertices we also use at least one arc $(k,j), j \neq i$ out of $k$.

In a similar fashion, artificial vertex reversal can occur at the end of a subaisle. For example, suppose in Figure \ref{fig4}  we have a small subtour for some trolley $t$ in the second aisle involving (10,22) and (22,10), so the trolley goes from vertex 10 to vertex 22 and back to vertex 10. This subtour clearly satisfies the cross-aisle constraint associated with the dotted cut shown, which here is $x_{t,21,4}+x_{t,22,10}+x_{t,23,16} \geq z_{ot}$. In order to prevent this situation we can impose the additional constraints $x_{t,22,21}+ x_{t,22,9}  + x_{t,22,23} \geq x_{t,10,22}$  (which ensures that if we use (10,22) we go onward from vertex 22 to appropriate adjacent vertices) and
$x_{t,21,22}+ x_{t,9,22}  + x_{t,23,22} \geq x_{t,22,10}$   (which ensures that if we use (22,10) we have come to  vertex 22 from appropriate adjacent vertices).
Generalising, we have that:

\begin{optprog}
& \sum_{\substack{j \in \delta^+(v(a,c)) \\ j \neq S(a,c,1) }} \kern -1.2em x_{t,v(a,c),j} & \;\geq\; & x_{t,S(a,c,1),v(a,c)} 
& \;\;\; \forall a \in \mathbb{A}, c \in \mathbb{C}\setminus[W_C], t \in {\cal T}
\label{jeb19} \\

& \sum_{\substack{j \in \delta^-(v(a,c)) \\ j \neq S(a,c,1) }} \kern -1.2em x_{t,j,v(a,c)} & \;\geq\; & x_{t,v(a,c),S(a,c,1)} 
&\;\;\; \forall a \in \mathbb{A}, c \in \mathbb{C}\setminus[W_C], t \in {\cal T}
\label{jeb20}\\

& \sum_{\substack{j \in \delta^+(v(a,c)) \\ j \neq S(a,c-1,R(a,c-1)) }} \kern -2.7em x_{t,v(a,c),j} & \;\geq\; & x_{t,S(a,c-1,R(a,c-1)),v(a,c)} 
&\;\;\; \forall a \in \mathbb{A}, c \in \mathbb{C}\setminus[1], t \in {\cal T}
\label{jeb21} 
\end{optprog}

\begin{optprog}
& \sum_{\substack{j \in \delta^-(v(a,c)) \\ j \neq S(a,c-1,R(a,c-1)) }} \kern -2.7em x_{t,j,v(a,c)} & \;\geq\; & x_{t,v(a,c),S(a,c-1,R(a,c-1))} 
&\;\;\; \forall a \in \mathbb{A}, c \in \mathbb{C}\setminus[1], t \in {\cal T}
\label{jeb22} 
\end{optprog}

Here equations (\ref{jeb19},\ref{jeb20}) deal with reversals with respect to an artificial vertex $v(a,c)$ and the adjacent southern  subaisle vertex $S(a,c,1)$.  
Equations (\ref{jeb21},\ref{jeb22}) deal with reversals with respect to an artificial vertex $v(a,c)$ and the adjacent northern  subaisle vertex $S(a,c-1,R(a,c-1))$.


\subsection{Pass through}

Suppose that order $o$ is not going to be picked by trolley $t$, so $z_{ot}=0$. Then if trolley $t$ does visit any vertices that 
\textbf{\emph{only}} exist in the reduced graph to service order $o$ it must pass through them (either from north to south or from south to north). This is because reversal at a vertex where no picking occurs means that we are travelling extra distance for no practical effect
(as for reversal at artificial vertices considered above). 
For example in Figure~\ref{fig4} suppose that trolley $t$ is not going to service the single order associated with vertex 5 in the circled subaisle. Then it must be true that $x_{t,4,5}=x_{t,5,6}$ and $x_{t,6,5}=x_{t,5,4}$. Here, for example, the first equality ensures that arcs (4,5) and (5,6) have equal value for trolley $t$, whether the trolley visits vertex 5 or not.

Generalising let $\Psi(S(a,c,r))$ be the set of orders associated with the subaisle vertex $S(a,c,r)$, the $r$'th vertex in the subaisle below artificial vertex $v(a,c)$. Then to ensure that we pass through a vertex at which no orders are picked we  have:

\begin{optprog}
& - \kern -1.5em \sum_{o \in \Psi(S(a,c,r))} \kern -1.5em z_{ot} & \;\leq\; & x_{t,S(a,c,r-1),S(a,c,r)} - x_{t,S(a,c,r),S(a,c,r+1)} & \;\leq\; \kern -1.5em \sum_{o \in \Psi(S(a,c,r))} \kern -1.5em z_{ot}, \label{jeb23} \\
& & & & r=1,\ldots,R(a,c), \forall a \in \mathbb{A}, c \in \mathbb{C}, t \in {\cal T} \nonumber \\

& - \kern -1.5em \sum_{o \in \Psi(S(a,c,r))} \kern -1.5em z_{ot} & \;\leq\; & x_{t,S(a,c,r+1),S(a,c,r)} - x_{t,S(a,c,r),S(a,c,r-1)} & \;\leq\; \kern -1.5em \sum_{o \in \Psi(S(a,c,r))} \kern -1.5em z_{ot}, \label{jeb24}\\
& & & & r=1,\ldots,R(a,c), \forall a \in \mathbb{A}, c \in \mathbb{C}, t \in {\cal T} \nonumber  
\end{optprog}





\subsection{Remark}

In this section we have introduced a significant number of valid inequalities (cuts) for the JOBPRP. To the best of our knowledge, none of these cuts has been presented before in the literature. In particular, none were introduced in our previous work \citep{valle2016}, where we focussed on different formulations instead of enhancing those formulations via the use of inequalities.

The number of additional constraints and additional variables needed for all valid inequalities are polynomially bounded. Computational results, reported in the next section, show that the addition of these inequalities greatly reduces the number of disconnected solutions found during the branch-and-cut procedure. This allowed us to switch off the search for violated connectivity constraints when connected fractional solutions are found, which requires the use of the (relatively) computationally expensive max-flow algorithm. Rather we only search for violated connectivity constraints when we encounter an all-integer solution.


\section{Computational results}
\label{sec:results}

We conducted a series of experiments to empirically evaluate the impact of adding the proposed valid inequalities defined in Section~\ref{sec:ineq} to the formulation introduced in Section~\ref{sec:IP}. The generation of test problems is explained below.

\subsection{Test problems}
\label{sec:test}

In order to simulate a realistic supermarket environment, we make use of the publicly available \cite{foodmart}. The database is composed of  anonymised customer purchases over two years for a chain of supermarkets. There are a total of $1560$ distinct products, separated in product classes containing 4 different category levels. It also contains approximately 270000 orders for the period 1997-1998, each composed of a customer id, a list of distinct products purchased, the number of items for each distinct product and the purchase date. 

\subsubsection{Warehouse}

No information about warehouse layouts and product placement exists in Foodmart, therefore we built a warehouse layout generator to simulate both. The generator is based on a previous version developed by \cite{theys2010} and kindly provided by Dr. Birger Raa. The generator creates warehouses that must be able to hold a minimum predetermined number of distinct products given a (fixed) number of aisles, cross-aisles and shelves. Arbitrary lengths are also given, in metres, for aisles and cross-aisles width, as well as rack depth and slot width. The distance from the origin to its closest artificial vertex (the black vertex in the top left corner of Figure \ref{fig2}) is also given.

The generator computes the number of slots a shelf must have in order to hold at least the required number of products, while keeping the number of empty slots to a minimum. As an example, the warehouse in Figure \ref{fig1} was computed for a minimum of 104 products, while having 3 aisles, 3 cross-aisles and 2 shelves. Each shelf must have at least 9 slots, so that the total number of individual products in the warehouse is 108 (four slots would be empty in this example). The generator also computes the position of cross-aisles such that aisles are divided in subaisles as equal (in terms of number of slots) as possible. The placement of products in slots is done by sorting all products from the highest category level to the lowest, and placing them in consecutive slots, so that similar products are close to each other. 

A single warehouse layout is used for the majority of the work reported here. It contains 8 aisles, 3 cross-aisles and 3 shelves. Each shelf holds 33 slots, so this warehouse can store 1584 distinct products (enough for all 1560 products in the Foodmart database). The distance from the origin to the first artificial vertex is 4m, the aisle and cross-aisle widths are 3m, and both the slot width and rack depth are 1m.

\subsubsection{Orders}

We observed in the Foodmart database that orders are generally very small (the vast majority containing up to only 4 or 5 distinct products). On the other hand, online orders (which inspired this problem) may be composed of dozens of items. 

To produce larger orders, we combined different Foodmart orders into a single one. For every customer, all of their purchases made in the first $\Delta$ days are combined into a single order. The combined order may contain not only more distinct products, but also a higher quantity of items of a single product. 

A test instance is taken as the $O$ orders with the highest number of distinct products. If $O = 5$, the 5 largest combined orders make up the test set; if $O = 6$, we take the same orders as in $O = 5$ plus the sixth largest combined order. We created several test instances for $\Delta = \{5, 10, 20\}$ and $O = \{5, \ldots, 30\}$. Table \ref{table1} shows the number of distinct products, the total number of items and the number of required baskets for the largest 30 orders (for each value of $\Delta$). We consider that each basket may carry up to $40$ items. All test instances are publicly available as mentioned previously above.

\subsubsection{Capacity of baskets and number of trolleys} With regard to parameter values we set $B = 8$, each basket holding a maximum of $40$ items, irrespective of their sizes or weights. For every test instance, we define the number of available trolleys $T = \Big\lceil \frac{\sum_{o \in O} b_o}{B} + 0.2 \Big\rceil $. Finding the exact minimum $T$ required to service all orders is an optimisation problem on its own, we however do not tackle this problem in this work. Not all trolleys need to be used as the solution may leave some trolleys idle. 

\begin{table}[!htbp]
\centering
{\scriptsize
\renewcommand{\tabcolsep}{1mm} 
\renewcommand{\arraystretch}{1} 
\begin{tabular}{|c|ccc|ccc|ccc|}
\hline
\multirow{2}{*}{Order Index} & \multicolumn{3}{c|}{$\Delta = 5$} & \multicolumn{3}{c|}{$\Delta = 10$} & \multicolumn{3}{c|}{$\Delta = 20$} \\
 & Products & Items & Baskets & Products & Items & Baskets & Products & Items & Baskets\\
\hline
   1 &  18  &    63   &  2 &  23   &   79  &  2  &  23   &   79  &  2 \\
   2 &  11  &    39   &  1 &  23   &   66  &  2  &  23   &   66  &  2 \\
   3 &  11  &    39   &  1 &  12   &   34  &  1  &  16   &   45  &  2 \\
   4 &  11  &    35   &  1 &  11   &   40  &  1  &  15   &   42  &  2 \\
   5 &   8  &    28   &  1 &  11   &   39  &  1  &  14   &   50  &  2 \\
   6 &   8  &    28   &  1 &  11   &   39  &  1  &  14   &   49  &  2 \\
   7 &   7  &    28   &  1 &  11   &   36  &  1  &  14   &   40  &  1 \\
   8 &   7  &    27   &  1 &  11   &   35  &  1  &  13   &   42  &  2 \\
   9 &   7  &    27   &  1 &  11   &   31  &  1  &  13   &   41  &  2 \\
  10 &   7  &    24   &  1 &  10   &   38  &  1  &  13   &   40  &  1 \\
  11 &   7  &    24   &  1 &  10   &   36  &  1  &  13   &   39  &  1 \\
  12 &   7  &    23   &  1 &  10   &   32  &  1  &  12   &   48  &  2 \\
  13 &   7  &    23   &  1 &  10   &   30  &  1  &  12   &   45  &  2 \\
  14 &   7  &    22   &  1 &  10   &   29  &  1  &  12   &   43  &  2 \\
  15 &   7  &    21   &  1 &   9   &   28  &  1  &  12   &   41  &  2 \\
  16 &   7  &    20   &  1 &   9   &   26  &  1  &  12   &   39  &  1 \\
  17 &   7  &    20   &  1 &   8   &   27  &  1  &  12   &   36  &  1 \\
  18 &   7  &    20   &  1 &   8   &   25  &  1  &  12   &   34  &  1 \\
  19 &   7  &    20   &  1 &   8   &   22  &  1  &  12   &   34  &  1 \\
  20 &   7  &    19   &  1 &   7   &   26  &  1  &  12   &   33  &  1 \\
  21 &   7  &    19   &  1 &   7   &   23  &  1  &  11   &   43  &  2 \\
  22 &   7  &    19   &  1 &   7   &   23  &  1  &  11   &   40  &  1 \\
  23 &   7  &    18   &  1 &   7   &   22  &  1  &  11   &   39  &  1 \\
  24 &   7  &    18   &  1 &   7   &   20  &  1  &  11   &   39  &  1 \\
  25 &   7  &    10   &  1 &   7   &   20  &  1  &  11   &   38  &  1 \\
  26 &   6  &    21   &  1 &   7   &   20  &  1  &  11   &   38  &  1 \\
  27 &   6  &    19   &  1 &   7   &   19  &  1  &  11   &   37  &  1 \\
  28 &   6  &    17   &  1 &   7   &   18  &  1  &  11   &   36  &  1 \\
  29 &   6  &    14   &  1 &   7   &   17  &  1  &  11   &   33  &  1 \\
  30 &   6  &     9   &  1 &   7   &   10  &  1  &  11   &   31  &  1 \\
\hline
\end{tabular}
\caption{Largest combined orders for $\Delta = \{5, 10, 20\}$}
\label{table1}
}
\end{table}

Details about test instances are shown in Table \ref{table2}. For each value of $\Delta = \{5, 10, 20\}$, there are five columns: The total number of baskets needed to carry all orders (labelled as $\sum_{o \in O} b_o$), the number of available trolleys $T$, the actual number of trolleys $T^*$ required to service all orders and the number of vertices $|V|$ and edges $|A|$ of each respective graph $D$. For higher values of $\Delta$, both the total number of items and the number of distinct products increase, this is reflected in the higher number of baskets required and the larger number of vertices and arcs for higher values of $\Delta$.

\begin{table}[!htbp]
\centering
{\scriptsize
\renewcommand{\tabcolsep}{1mm} 
\renewcommand{\arraystretch}{1.4} 
\begin{tabular}{|c|ccccc|ccccc|ccccc|}
\hline
\multirow{2}{*}{$O$} & \multicolumn{5}{c|}{$\Delta = 5$} & \multicolumn{5}{c|}{$\Delta = 10$} & \multicolumn{5}{c|}{$\Delta = 20$} \\
 & $\sum_{o \in O} b_o$ & $T$ & $T^*$ & $|V|$ & $|A|$ & $\sum_{o \in O} b_o$ & $T$ & $T^*$ & $|V|$ & $|A|$ & $\sum_{o \in O} b_o$ & $T$ & $T^*$ & $|V|$ & $|A|$ \\
\hline
    5 &     6 &     1 &     1 &    76 &   192 &     7 &     2 &     1 &    90 &   220 &    10 &     2 &     2 &   104 &   248 \\
    6 &     7 &     2 &     1 &    82 &   204 &     8 &     2 &     1 &   101 &   242 &    11 &     2 &     2 &   114 &   268 \\
    7 &     8 &     2 &     1 &    89 &   218 &     9 &     2 &     2 &   109 &   258 &    13 &     2 &     2 &   123 &   286 \\
    8 &     9 &     2 &     2 &    95 &   230 &    10 &     2 &     2 &   114 &   268 &    14 &     2 &     2 &   133 &   306 \\
    9 &    10 &     2 &     2 &   101 &   242 &    11 &     2 &     2 &   123 &   286 &    16 &     3 &     2 &   141 &   322 \\
   10 &    11 &     2 &     2 &   106 &   252 &    12 &     2 &     2 &   131 &   302 &    17 &     3 &     3 &   147 &   334 \\
   11 &    12 &     2 &     2 &   111 &   262 &    13 &     2 &     2 &   140 &   320 &    19 &     3 &     3 &   149 &   338 \\
   12 &    13 &     2 &     2 &   116 &   272 &    14 &     2 &     2 &   145 &   330 &    20 &     3 &     3 &   153 &   346 \\
   13 &    14 &     2 &     2 &   121 &   282 &    15 &     3 &     2 &   150 &   340 &    21 &     3 &     3 &   157 &   354 \\
   14 &    15 &     3 &     2 &   128 &   296 &    16 &     3 &     2 &   156 &   352 &    23 &     4 &     3 &   162 &   364 \\
   15 &    16 &     3 &     2 &   132 &   304 &    17 &     3 &     3 &   160 &   360 &    24 &     4 &     3 &   169 &   378 \\
   20 &    21 &     3 &     3 &   153 &   346 &    22 &     3 &     3 &   179 &   398 &    32 &     5 &     4 &   196 &   432 \\
   25 &    26 &     4 &     4 &   174 &   388 &    27 &     4 &     4 &   195 &   430 &    38 &     5 &     5 &   211 &   462 \\
   30 &    31 &     5 &     4 &   186 &   412 &    32 &     5 &     4 &   207 &   454 &    43 &     6 &     6 &   228 &   496 \\
\hline
\end{tabular}
\caption{Instances features}
\label{table2}
}
\end{table}

\subsection{Addition of all valid inequalities}

In this section, we present computational results to analyse the impact of adding all the valid inequalities introduced in Section~\ref{sec:ineq}. We use {\it\textbf{original model}} to refer to formulation (\ref{eq1})-(\ref{eq17}) and {\it\textbf{reinforced model}} to refer to the same formulation with the inclusion of inequalities (\ref{jeb1})-(\ref{jeb9}) and (\ref{jeb11})-(\ref{jeb24}). The original model is solved via the FBC algorithm and the reinforced model is solved via the IBC algorithm - as defined in Section~\ref{sec:branch}.

We used an Intel Xeon \@ 2.40GHz with 32GB of RAM and Linux as the operating system. The code was written in C++ and CPLEX 12.6.0 \citep{cplex} was used as the mixed-integer solver. The heuristic described in Section~\ref{sec:impl} is used to warm-start the algorithm with valid JOBPRP upper bounds and solver parameters are set as described in the same section. A maximum time limit of 6 CPU hours was imposed to each algorithm and instance.

Table \ref{table3} shows results for the instances described in Table \ref{table2}. For both the reinforced and original models, we include seven columns. {\it\textbf{T(s)}} denotes the total elapsed time in seconds, {\it\textbf{UB}} and {\it\textbf{LB}} respectively represent the best upper and lower bounds obtained at the end of the search either when the instance was solved to proven optimality or when the time limit was reached. {\it\textbf{GAP}} is defined as $100 (\text{UB} - \text{LB})/\text{UB}$. {\it\textbf{FLB}} is the lower bound obtained at the end of the root node of the branch-and-bound search tree and {\it\textbf{FGAP}} is defined as $100 (\text{UB} - \text{FLB})/\text{UB}$. That is, for instances solved to optimality, {\it\textbf{FGAP}} denotes the gap between the optimal solution and the lower bound obtained after solving the root node, which gives an indication of how strong the linear relaxation is. Finally, {\it\textbf{NS}} stands for the total number of nodes investigated during the search. A bold value in the table indicates a better value in comparison to the other corresponding column.

\begin{table}[!htbp]
\centering
{\scriptsize
\renewcommand{\tabcolsep}{1mm} 
\renewcommand{\arraystretch}{1.4} 
\begin{tabular}{|cc|rrrrrrr|rrrrrrr|}
\hline
\multirow{2}{*}{$\Delta$} & \multirow{2}{*}{$O$} & \multicolumn{7}{c|}{Original model solved with FBC} & \multicolumn{7}{c|}{Reinforced model solved with IBC}\\
 &  & \multicolumn{1}{c}{T(s)} & \multicolumn{1}{c}{UB} & \multicolumn{1}{c}{GAP} & \multicolumn{1}{c}{LB} & \multicolumn{1}{c}{FGAP} & \multicolumn{1}{c}{FLB} & \multicolumn{1}{c|}{NS} & \multicolumn{1}{c}{T(s)} & \multicolumn{1}{c}{UB} & \multicolumn{1}{c}{GAP} & \multicolumn{1}{c}{LB} & \multicolumn{1}{c}{FGAP} & \multicolumn{1}{c}{FLB} & \multicolumn{1}{c|}{NS}\\
\hline
   5 &    5 & \textbf{    2.6} &   348.6 & \multicolumn{1}{c}{--} & \multicolumn{1}{c}{--} &  0.6 &    346.5 &     7 &    3.9 &   348.6 & \multicolumn{1}{c}{--} & \multicolumn{1}{c}{--} &\textbf{  0.6} &    346.5 &    36 \\
 &    6 &    14.4 &   364.8 & \multicolumn{1}{c}{--} & \multicolumn{1}{c}{--} &  0.0 &    364.8 &     1 &\textbf{    2.9} &   364.8 & \multicolumn{1}{c}{--} & \multicolumn{1}{c}{--} &  0.0 &    364.8 &     1 \\
 &    7 &    13.6 &   374.8 & \multicolumn{1}{c}{--} & \multicolumn{1}{c}{--} &\textbf{  0.5} &    372.9 &    45 &\textbf{    8.9} &   374.8 & \multicolumn{1}{c}{--} & \multicolumn{1}{c}{--} &  0.5 &    372.8 &   141 \\
 &    8 &   149.0 &   503.8 & \multicolumn{1}{c}{--} & \multicolumn{1}{c}{--} & 19.2 &    407.1 &   559 &\textbf{   95.2} &   503.8 & \multicolumn{1}{c}{--} & \multicolumn{1}{c}{--} &\textbf{  8.3} &    462.1 &   146 \\
 &    9 &   291.2 &   539.6 & \multicolumn{1}{c}{--} & \multicolumn{1}{c}{--} & 22.9 &    416.1 &  2622 &\textbf{  151.8} &   539.6 & \multicolumn{1}{c}{--} & \multicolumn{1}{c}{--} &\textbf{ 11.2} &    479.2 &   268 \\
 &   10 &  1072.3 &   581.4 & \multicolumn{1}{c}{--} & \multicolumn{1}{c}{--} & 26.2 &    429.3 & 13445 &\textbf{  111.0} &   581.4 & \multicolumn{1}{c}{--} & \multicolumn{1}{c}{--} &\textbf{ 14.4} &    497.8 &    88 \\
 &   11 &  5064.9 &   613.5 & \multicolumn{1}{c}{--} & \multicolumn{1}{c}{--} & 28.4 &    439.3 & 77291 &\textbf{   97.1} &   613.5 & \multicolumn{1}{c}{--} & \multicolumn{1}{c}{--} &\textbf{ 16.9} &    509.9 &  1009 \\
 &   12 &  6129.9 &   621.4 & \multicolumn{1}{c}{--} & \multicolumn{1}{c}{--} & 26.4 &    457.1 & 84901 &\textbf{  256.7} &   621.4 & \multicolumn{1}{c}{--} & \multicolumn{1}{c}{--} &\textbf{ 15.0} &    527.9 &  2261 \\
 &   13 &  3613.5 &   623.4 & \multicolumn{1}{c}{--} & \multicolumn{1}{c}{--} & 25.4 &    465.1 & 42659 &\textbf{  168.9} &   623.4 & \multicolumn{1}{c}{--} & \multicolumn{1}{c}{--} &\textbf{ 15.3} &    528.2 &   395 \\
 &   14 & \multicolumn{1}{c}{--} &   639.4 &   2.9 &    620.6 &  26.6 &    469.1 & 62700 &\textbf{  263.8} &   639.3 & \multicolumn{1}{c}{--} & \multicolumn{1}{c}{--} &\textbf{ 17.8} &    525.5 &   995 \\
 &   15 & \multicolumn{1}{c}{--} &   667.8 &   8.6 &    610.5 &  28.5 &    477.2 & 59000 &\textbf{  348.9} &   653.4 & \multicolumn{1}{c}{--} & \multicolumn{1}{c}{--} &\textbf{ 20.0} &    522.4 &   975 \\
 &   20 & \multicolumn{1}{c}{--} &   946.4 &  30.7 &    655.6 &  45.4 &    516.8 & 38500 &\textbf{ 2990.9} &   870.4 & \multicolumn{1}{c}{--} & \multicolumn{1}{c}{--} &\textbf{ 31.2} &    598.4 & 28931 \\
 &   25 & \multicolumn{1}{c}{--} &  1248.4 &  48.2 &    646.7 &  56.0 &    549.0 & 15800 &\multicolumn{1}{c}{--} &\textbf{  1127.4} & \textbf{ 17.8} &    927.0 & \textbf{ 43.2} &    640.8 & 123900 \\
 &   30 & \multicolumn{1}{c}{--} &  1509.3 &  58.3 &    628.7 &  63.2 &    555.4 &  8963 &\multicolumn{1}{c}{--} &\textbf{  1221.9} & \textbf{ 26.8} &    894.3 & \textbf{ 47.0} &    647.5 & 70400 \\
\multicolumn{2}{|r|}{\textbf{Avg}} & \textbf{8882.2} &  & \textbf{10.6} &  & \textbf{26.4} &  \multicolumn{2}{r|}{\textbf{29035.2}} & \textbf{3407.1} &  & \textbf{3.2} &  & \textbf{17.2} &  \multicolumn{2}{r|}{\textbf{16396.1}}\\
\cline{1-16}
  10 &    5 &     4.5 &   371.1 & \multicolumn{1}{c}{--} & \multicolumn{1}{c}{--} &  0.0 &    371.1 &     1 &\textbf{    1.5} &   371.1 & \multicolumn{1}{c}{--} & \multicolumn{1}{c}{--} &  0.0 &    371.1 &     1 \\
 &    6 &    16.4 &   377.1 & \multicolumn{1}{c}{--} & \multicolumn{1}{c}{--} &\textbf{  0.6} &    374.7 &   123 &\textbf{    6.6} &   377.1 & \multicolumn{1}{c}{--} & \multicolumn{1}{c}{--} &  0.6 &    374.7 &   207 \\
 &    7 & \textbf{  190.4} &   549.8 & \multicolumn{1}{c}{--} & \multicolumn{1}{c}{--} & 22.4 &    426.8 &   848 &  197.5 &   549.8 & \multicolumn{1}{c}{--} & \multicolumn{1}{c}{--} &\textbf{  8.4} &    503.4 &    13 \\
 &    8 &   331.1 &   584.2 & \multicolumn{1}{c}{--} & \multicolumn{1}{c}{--} & 24.6 &    440.4 &  2222 &\textbf{  286.8} &   584.2 & \multicolumn{1}{c}{--} & \multicolumn{1}{c}{--} &\textbf{ 13.0} &    508.4 &   110 \\
 &    9 &  5529.7 &   637.4 & \multicolumn{1}{c}{--} & \multicolumn{1}{c}{--} & 28.6 &    455.1 & 96788 &\textbf{  356.1} &   637.4 & \multicolumn{1}{c}{--} & \multicolumn{1}{c}{--} &\textbf{ 17.1} &    528.6 &  8608 \\
 &   10 & 10512.3 &   661.8 & \multicolumn{1}{c}{--} & \multicolumn{1}{c}{--} & 29.0 &    469.9 & 136971 &\textbf{  393.0} &   661.8 & \multicolumn{1}{c}{--} & \multicolumn{1}{c}{--} &\textbf{ 16.4} &    553.4 & 10439 \\
 &   11 & \multicolumn{1}{c}{--} &   699.8 &   5.4 &    662.3 &  31.1 &    482.3 & 140888 &\textbf{  383.7} &   699.8 & \multicolumn{1}{c}{--} & \multicolumn{1}{c}{--} &\textbf{ 20.9} &    553.8 & 21044 \\
 &   12 & \multicolumn{1}{c}{--} &   719.8 &   6.6 &    672.0 &  31.4 &    493.6 & 110415 &\textbf{  217.7} &   707.7 & \multicolumn{1}{c}{--} & \multicolumn{1}{c}{--} &\textbf{ 21.8} &    553.6 &  2481 \\
 &   13 & \multicolumn{1}{c}{--} &   733.8 &   8.7 &    670.1 &  31.3 &    504.1 & 74400 &\textbf{  385.4} &   725.7 & \multicolumn{1}{c}{--} & \multicolumn{1}{c}{--} &\textbf{ 22.8} &    560.2 & 10370 \\
 &   14 & \multicolumn{1}{c}{--} &   739.8 &   8.7 &    675.5 &  31.2 &    509.1 & 69600 &\textbf{  483.8} &   727.8 & \multicolumn{1}{c}{--} & \multicolumn{1}{c}{--} &\textbf{ 22.7} &    562.4 &  8597 \\
 &   15 & \multicolumn{1}{c}{--} &   933.0 &  28.1 &    670.8 &  43.5 &    527.2 & 39600 &\textbf{ 1832.6} &   882.6 & \multicolumn{1}{c}{--} & \multicolumn{1}{c}{--} &\textbf{ 32.8} &    593.1 & 26398 \\
 &   20 & \multicolumn{1}{c}{--} &  1093.0 &  34.5 &    716.1 &  48.9 &    559.0 & 33900 &\textbf{11015.4} &   992.4 & \multicolumn{1}{c}{--} & \multicolumn{1}{c}{--} &\textbf{ 32.3} &    672.0 & 70496 \\
 &   25 & \multicolumn{1}{c}{--} &  1300.4 &  44.6 &    720.0 &  54.2 &    595.4 & 14219 &\multicolumn{1}{c}{--} &\textbf{  1213.4} & \textbf{ 16.7} &   1011.3 & \textbf{ 41.7} &    707.8 & 133600 \\
 &   30 & \multicolumn{1}{c}{--} &  1565.7 &  55.8 &    692.1 &  60.9 &    612.8 &  7439 &\multicolumn{1}{c}{--} &\textbf{  1330.0} & \textbf{ 27.6} &    963.1 & \textbf{ 46.9} &    706.3 & 74100 \\
\multicolumn{2}{|r|}{\textbf{Avg}} & \textbf{13527.5} &  & \textbf{13.7} &  & \textbf{31.3} &  \multicolumn{2}{r|}{\textbf{51958.1}} & \textbf{4197.2} &  & \textbf{3.2} &  & \textbf{21.2} &  \multicolumn{2}{r|}{\textbf{26176.0}}\\
\cline{1-16}
  20 &    5 &   212.8 &   573.8 & \multicolumn{1}{c}{--} & \multicolumn{1}{c}{--} & 20.1 &    458.4 &  1325 &\textbf{  110.9} &   573.8 & \multicolumn{1}{c}{--} & \multicolumn{1}{c}{--} &\textbf{  6.7} &    535.3 &   233 \\
 &    6 &   967.7 &   656.2 & \multicolumn{1}{c}{--} & \multicolumn{1}{c}{--} & 27.1 &    478.4 & 14378 &\textbf{  278.9} &   656.2 & \multicolumn{1}{c}{--} & \multicolumn{1}{c}{--} &\textbf{ 11.5} &    580.9 &  1164 \\
 &    7 &  5123.9 &   689.8 & \multicolumn{1}{c}{--} & \multicolumn{1}{c}{--} & 27.4 &    501.1 & 83100 &\textbf{  263.2} &   689.8 & \multicolumn{1}{c}{--} & \multicolumn{1}{c}{--} &\textbf{ 17.2} &    571.4 &   579 \\
 &    8 & 10929.6 &   697.8 & \multicolumn{1}{c}{--} & \multicolumn{1}{c}{--} & 28.2 &    501.2 & 131230 &\textbf{  189.0} &   697.8 & \multicolumn{1}{c}{--} & \multicolumn{1}{c}{--} &\textbf{ 18.5} &    568.7 &  2886 \\
 &    9 & \multicolumn{1}{c}{--} &   727.7 &   3.2 &    704.5 &  29.6 &    512.3 & 112033 &\textbf{  305.5} &   727.7 & \multicolumn{1}{c}{--} & \multicolumn{1}{c}{--} &\textbf{ 20.6} &    578.0 &  3620 \\
 &   10 & \multicolumn{1}{c}{--} &   935.3 &  16.1 &    784.2 &  43.4 &    529.5 & 62300 &\textbf{  584.2} &   920.5 & \multicolumn{1}{c}{--} & \multicolumn{1}{c}{--} &\textbf{ 31.1} &    634.2 &  7198 \\
 &   11 & \multicolumn{1}{c}{--} &  1005.2 &  21.1 &    792.8 &  45.4 &    548.5 & 53700 &\textbf{  703.8} &   980.5 & \multicolumn{1}{c}{--} & \multicolumn{1}{c}{--} &\textbf{ 31.4} &    672.9 & 10620 \\
 &   12 & \multicolumn{1}{c}{--} &  1047.3 &  24.2 &    793.4 &  47.6 &    548.5 & 54000 &\textbf{  830.1} &  1004.3 & \multicolumn{1}{c}{--} & \multicolumn{1}{c}{--} &\textbf{ 33.7} &    666.0 & 20221 \\
 &   13 & \multicolumn{1}{c}{--} &  1048.6 &  23.7 &    799.8 &  47.2 &    553.2 & 59400 &\textbf{ 1177.8} &  1009.1 & \multicolumn{1}{c}{--} & \multicolumn{1}{c}{--} &\textbf{ 31.0} &    696.5 & 30072 \\
 &   14 & \multicolumn{1}{c}{--} &  1079.3 &  28.4 &    772.9 &  47.8 &    563.1 & 31400 &\textbf{ 1306.0} &  1011.1 & \multicolumn{1}{c}{--} & \multicolumn{1}{c}{--} &\textbf{ 32.8} &    679.2 &  8097 \\
 &   15 & \multicolumn{1}{c}{--} &  1135.3 &  34.0 &    749.4 &  50.0 &    567.3 & 25292 &\textbf{ 3793.7} &  1028.7 & \multicolumn{1}{c}{--} & \multicolumn{1}{c}{--} &\textbf{ 34.7} &    671.5 & 38243 \\
 &   20 & \multicolumn{1}{c}{--} &  1428.4 &  47.3 &    753.3 &  57.1 &    612.9 &  8560 &\multicolumn{1}{c}{--} &\textbf{  1335.2} & \textbf{ 15.8} &   1124.6 & \textbf{ 45.4} &    729.0 & 67600 \\
 &   25 & \multicolumn{1}{c}{--} &  1777.5 &  57.2 &    760.1 &  64.0 &    639.6 &  8596 &\multicolumn{1}{c}{--} &\textbf{  1694.9} & \textbf{ 31.2} &   1166.8 & \textbf{ 51.4} &    824.5 & 67200 \\
 &   30 & \multicolumn{1}{c}{--} &  2076.0 &  64.3 &    741.6 &  68.4 &    656.9 &  3225 &\multicolumn{1}{c}{--} &\textbf{  1966.9} & \textbf{ 45.5} &   1071.9 & \textbf{ 56.6} &    853.4 & 33400 \\
\multicolumn{2}{|r|}{\textbf{Avg}} & \textbf{16659.6} &  & \textbf{22.8} &  & \textbf{43.1} &  \multicolumn{2}{r|}{\textbf{46324.2}} & \textbf{5310.2} &  & \textbf{6.6} &  & \textbf{30.2} &  \multicolumn{2}{r|}{\textbf{20795.2}}\\
\cline{1-16}
\hline
\multicolumn{2}{r}{\textbf{Avg all}} & \textbf{13023.1} &  & \textbf{15.7} &  & \textbf{33.6} &  \multicolumn{2}{r}{\textbf{42439.2}} & \textbf{4304.8} &  & \textbf{4.3} &  & \textbf{22.9} &  \multicolumn{2}{r}{\textbf{21122.5}}\end{tabular}
}
\caption{Solving the original and reinforced models using FBC and IBC respectively. A symbol ``$-$'' for \textbf{T(s)} entries indicates that the instance was not solved to proven optimality within the time limit.}
\label{table3}
\end{table}

For each value of $\Delta$, we present a further row with the average values for T(s), GAP, FGAP and NS. The last row in the table shows an ensemble average for all instances together. 

The FBC results for the original model are slightly improved when compared to the results in \cite{valle2016}. This is partly due to the addition of the $\alpha$ variables, which also made possible the introduction of symmetry-breaking constraints (\ref{eq17}), and partly due to the solver presolve being enabled - as opposed to our previous work. 

It is clear from the table that the addition of valid inequalities proposed in Section~\ref{sec:ineq} greatly improves overall performance. The FBC based on the original model was able to solve 19 of the 42 instances to optimality within 6 hours, while the IBC for the reinforced model solved 35 out of 42 instances. As an example, for $\Delta = 5$, $O = 20$, the original model obtained a final optimality gap of $30.7\%$ after 6 hours, while the reinforced model solved it in about 50 minutes. Similar cases can be seen for other instances. 

An important difference between the two algorithms in Table \ref{table3} is the choice of FBC for the original model and IBC for the reinforced model. The push-relabel method by \cite{goldberg1988} for solving the max-flow problem has a relatively high computational cost and is thus responsible for a significant proportion of the time spent by the FBC. For instance, for the case where $\Delta = 20$ and $O = 30$, about 55\% of the 6 hours were spent in more than 2 million push-relabel executions. Nevertheless, we obtained worse results with the original model when using the IBC algorithm. These results are summarised in Table \ref{table4}, where due to space constraints we only include the rows with average values for each $\Delta$ and the ensemble average. 

\begin{table}[!htbp]
\centering
{\scriptsize
\renewcommand{\tabcolsep}{1.9mm} 
\renewcommand{\arraystretch}{1.5} 
\begin{tabular}{|c|rrrr|rrrr|}
\hline
\multirow{2}{*}{$\Delta$} & \multicolumn{4}{c|}{FBC} & \multicolumn{4}{c|}{IBC}\\
 & \multicolumn{1}{c}{T(s)} & \multicolumn{1}{c}{GAP} & \multicolumn{1}{c}{FGAP} & \multicolumn{1}{c|}{NS} & \multicolumn{1}{c}{T(s)} & \multicolumn{1}{c}{GAP} & \multicolumn{1}{c}{FGAP} & \multicolumn{1}{c|}{NS}\\
\hline
$5$  & 8877.2  &  10.6 &  26.4 &  29035.2 & 10460.1 & 13.9 &  38.9 & 815295.9  \\
$10$ & 13527.5 &  13.7 &  31.3 &  51958.1 & 14700.5 & 18.3 &  44.7 & 1574201.2 \\
$20$ & 16659.6 &  22.8 &  43.1 &  46324.2 & 18326.8 & 29.3 &  58.1 & 2284228.6 \\
\hline
\multicolumn{1}{c}{\textbf{Avg}} & \textbf{13021.4} & \textbf{15.7} & \textbf{33.6} & \multicolumn{1}{r}{\textbf{42439.2}} & \textbf{14495.8} & \textbf{20.5} & \textbf{47.2} & \multicolumn{1}{r}{\textbf{1557908.5}} \\

\end{tabular}
}
\caption{Comparison of the FBC and IBC algorithms for the original model.}
\label{table4}
\end{table}

Despite being time-consuming, as shown by comparing the NS columns for both FBC and IBC, the separation of constraints (\ref{eq11}) for the original model provides stronger linear relaxation bounds. This can be seen in Table \ref{table4} as both the GAP and FGAP values are on average lower for FBC than for IBC.

The situation is different for the reinforced model. Several of the inequalities included in the model naturally prevent many subtour solutions and the benefits of separating constraints (\ref{eq11}) for candidate fractional solutions are not enough to outweigh the computational cost. These results are summarised in Table \ref{table5}.

\begin{table}[!htbp]
\centering
{\scriptsize
\renewcommand{\tabcolsep}{1.9mm} 
\renewcommand{\arraystretch}{1.5} 
\begin{tabular}{|c|rrrr|rrrr|}
\hline
\multirow{2}{*}{$\Delta$} & \multicolumn{4}{c|}{FBC} & \multicolumn{4}{c|}{IBC}\\
 & \multicolumn{1}{c}{T(s)} & \multicolumn{1}{c}{GAP} & \multicolumn{1}{c}{FGAP} & \multicolumn{1}{c|}{NS} & \multicolumn{1}{c}{T(s)} & \multicolumn{1}{c}{GAP} & \multicolumn{1}{c}{FGAP} & \multicolumn{1}{c|}{NS}\\
\hline
$5$  & 4040.8  & 4.3  & 16.7 & 4069.8 & 3407.1  & 3.2  & 17.2 & 16396.1 \\
$10$ & 4756.4  & 4.4  & 21.5 & 7198.9 & 4197.2  & 3.2  & 21.2 & 26176.0 \\
$20$ & 6169.5  & 8.2  & 30.5 & 6027.1 & 5310.2  & 6.6  & 30.2 & 20795.2 \\
\hline
\multicolumn{1}{c}{\textbf{Avg}} & \textbf{4988.9} &  \textbf{5.6} & \textbf{22.9} & \multicolumn{1}{r}{\textbf{5765.3}} & \textbf{4304.8} & \textbf{4.3} &  \textbf{22.9} &  \multicolumn{1}{r}{\textbf{21122.5}} \\

\end{tabular}
}
\caption{Comparison of the FBC and IBC algorithms for the reinforced model.}
\label{table5}
\end{table}

While in Table \ref{table4} it is clear that FBC produces stronger linear relaxations (as the average FGAP is $33.6\%$ compared to $47.2\%$), the same is not true in Table \ref{table5} as the values for FGAP for both IBC and FBC are the same, at $22.9\%$. The average values for both T(s) and GAP are lower for IBC. Although the reinforced model still requires constraints (\ref{eq11}) to prevent disconnected solutions, it is not worth separating fractional ones as quite often no violated constraint is found.

It should also be noted that as some built-in CPLEX cuts were set to aggressive, the solver spends a reasonable amount of time in the solution of the root node. For the smaller instances, it is possible to solve them more quickly by setting these cuts to default. These instances are $O \leq 13$ for $\Delta = 5$, $O \leq 12$ for $\Delta = 10$ and $O \leq 8$ for $\Delta = 20$. For all other instances, the aggressive settings help reduce computation times (for the instances solved to optimality) and optimality gaps (for unsolved instances after 6 hours). 

\subsection{Comparison with other formulations}

In our previous work \citep{valle2016}, we introduced two other formulations for JOBPRP, both compact and based respectively on single-commodity flows and multi-commodity flows. The latter, being larger, was shown not to be competitive as the solver was not able to solve even the root node for large instances. However, in that work, the single-commodity flow formulation was empirically shown to be competitive when compared to the non-compact formulation. 

However, with the addition of inequalities (\ref{jeb1})-(\ref{jeb9}) and (\ref{jeb11})-(\ref{jeb24}), the non-compact formulation benefits from not having as many disconnected solutions, thus making the separation constraints (\ref{eq11}) for fractional solutions unnecessary. Due to being compact and not requiring an explicit branch-and-cut implementation, the same benefit is not seen for the single-commodity flow formulation, and thus we opted in this work to not include any compact formulation. Nevertheless, Table \ref{table6} shows summarised results comparing the reinforced versions of both formulations, under the exact same conditions. 

\begin{table}[!htbp]
\centering
{\scriptsize
\renewcommand{\tabcolsep}{1.9mm} 
\renewcommand{\arraystretch}{1.5} 
\begin{tabular}{|c|rrrr|rrrr|}
\hline
\multirow{2}{*}{$\Delta$} 
& \multicolumn{4}{>{\centering\arraybackslash}m{4.5cm}|}{Non-compact model (\ref{eq1})-(\ref{eq17})}
& \multicolumn{4}{>{\centering\arraybackslash}m{4.5cm}|}{Single-commodity flow model as in \citep{valle2016}}\\
 & \multicolumn{1}{c}{T(s)} & \multicolumn{1}{c}{GAP} & \multicolumn{1}{c}{FGAP} & \multicolumn{1}{c|}{NS} & \multicolumn{1}{c}{T(s)} & \multicolumn{1}{c}{GAP} & \multicolumn{1}{c}{FGAP} & \multicolumn{1}{c|}{NS}\\
\hline
$5$  & 3407.1  & 3.2  & 17.2 & 16396.1 & 5293.4  & 7.7  & 17.7 & 1988.6 \\
$10$ & 4197.2  & 3.2  & 21.2 & 26176.0 & 8077.8  & 7.9  & 22.2 & 10371.0\\
$20$ & 5310.2  & 6.6  & 30.2 & 20795.2 & 10130.4 & 12.2 & 31.4 & 10462.6\\
\hline
\multicolumn{1}{c}{\textbf{Avg}} & \textbf{4304.8} & \textbf{4.3} &  \textbf{22.9} &  \multicolumn{1}{r}{\textbf{21122.5}}& \textbf{7833.9} & \textbf{9.3} &  \textbf{23.8} &  \multicolumn{1}{r}{\textbf{7607.4}} \\

\end{tabular}
}
\caption{Comparison of the non-compact model, solved with IBC, and the single-commodity flow model as in \citep{valle2016}; both reinforced with (\ref{jeb1})-(\ref{jeb9}) and (\ref{jeb11})-(\ref{jeb24})}
\label{table6}
\end{table}

Both the average T(s) and GAP values are considerably lower for the non-compact formulation. Moreover, 4 instances previously solved to optimality were not solved in 6 hours with the single-commodity flow model: $O = 15$ for $\Delta = \{10, 20\}$ and $O = 20$ for $\Delta = \{5, 10\}$.

\subsection{Effect of individual cuts}

We have so far analysed the addition of all valid inequalities (cuts) defined in Section~\ref{sec:ineq}. This section aims to provide insight as to the individual contribution of each cut.

Table \ref{table7} shows summarised results for a subset of instances including the reinforced model, the original model and the reinforced model minus each one of the valid inequalities. Only instances that were solved to optimality by the reinforced model are included in the calculations shown in the table, namely $O \leq 20$ for $\Delta = \{5, 10\}$ and $O \leq 15$ for $\Delta = 20$. In the table, {\it \textbf{ FCAV}} stands for first cross-aisle vertex cuts (\ref{jeb2})-(\ref{jeb4}), {\it \textbf{ CA}} for cross-aisle cuts (\ref{jeb5})-(\ref{jeb8}), {\it \textbf{ A}} for aisle cuts (\ref{jeb9}), {\it \textbf{ SUB}} for subaisle cuts (\ref{jeb11})-(\ref{jeb17}), {\it \textbf{ AVR}} for artificial vertex reversal cuts (\ref{jeb18})-(\ref{jeb22}) and {\it \textbf{ PT}} for pass-through cuts (\ref{jeb23},\ref{jeb24}).

\begin{table}[!htbp]
\centering
{\scriptsize
\renewcommand{\tabcolsep}{1.5mm} 
\renewcommand{\arraystretch}{1.5} 
\begin{tabular}{|c|rrr|rrr|rrr|rrr|}
\hline
\multirow{2}{*}{Model} & \multicolumn{3}{c|}{$\Delta = 5, O \leq 20$} & \multicolumn{3}{c|}{$\Delta = 10, O \leq 20$} & \multicolumn{3}{c|}{$\Delta = 20, O \leq 15$} & \multicolumn{3}{c|}{Average}\\
 & \multicolumn{1}{c}{T(s)} & \multicolumn{1}{c}{GAP} & \multicolumn{1}{c|}{FGAP} & \multicolumn{1}{c}{T(s)} & \multicolumn{1}{c}{GAP} & \multicolumn{1}{c|}{FGAP} & \multicolumn{1}{c}{T(s)} & \multicolumn{1}{c}{GAP} & \multicolumn{1}{c|}{FGAP} & \multicolumn{1}{c}{T(s)} & \multicolumn{1}{c}{GAP} & \multicolumn{1}{c|}{FGAP} \\
\hline
\multicolumn{1}{|l|}{Reinforced model} 
   & 375.0  & 0.0  & 12.6 & 1296.7  & 0.0 & 17.4 & 867.6   & 0.0  & 24.5 & \textbf{845.8}   & \textbf{0.0} & \textbf{18.0} \\
\multicolumn{1}{|l|}{No FCAV cuts}
   & 596.9  & 0.0  & 13.9 & 680.7   & 0.0 & 17.9 & 987.5   & 0.0  & 25.2 & \textbf{748.4}   & \textbf{0.0} & \textbf{18.8} \\
\multicolumn{1}{|l|}{No CA cuts}
   & 342.7  & 0.0  & 12.4 & 888.3   & 0.0 & 17.9 & 788.3   & 0.0  & 24.6 & \textbf{669.8}   & \textbf{0.0} & \textbf{18.1} \\
\multicolumn{1}{|l|}{No A cuts}
   & 1954.8 & 0.1  & 17.0 & 7587.5  & 0.7 & 22.8 & 7752.1  & 0.2  & 29.8 & \textbf{5708.0}  & \textbf{0.3} & \textbf{23.0} \\
\multicolumn{1}{|l|}{No SUB cuts}
   & 2178.3 & 2.1  & 20.4 & 7370.2  & 3.3 & 29.0 & 10821.3 & 4.6  & 38.7 & \textbf{6674.7}  & \textbf{3.3} & \textbf{29.1} \\
\multicolumn{1}{|l|}{No AVR cuts}
   & 1269.6 & 0.0  & 14.3 & 1801.4  & 0.0 & 18.8 & 856.9   & 0.0  & 26.2 & \textbf{1322.2}  & \textbf{0.0} & \textbf{19.6} \\
\multicolumn{1}{|l|}{No PT cuts}
   & 388.8  & 0.0  & 13.2 & 827.1   & 0.0 & 18.3 & 1337.7  & 0.0  & 25.1 & \textbf{837.3}   & \textbf{0.0} & \textbf{18.7} \\
\multicolumn{1}{|l|}{Original model}
   & 6756.7 & 3.5  & 20.8 & 12182.0 & 7.7 & 26.9 & 15312.2 & 13.7 & 37.6 & \textbf{11305.7} & \textbf{8.1} & \textbf{28.2} \\
\hline
\end{tabular}
}
\caption{Comparison of the reinforced model, original model and the reinforced model minus each valid inequality separately. Only instances solved to optimality with the reinforced model are considered. All models containing valid inequalities were solved with IBC, the original model was solved with FBC.}
\label{table7}
\end{table}

The table suggests that the most effective cuts are SUB and A cuts. Their importance is highlighted by the significant increase in the average T(s) and FGAP values when we remove each of these separately. When removing A cuts only, 7 instances solved by the reinforced model could not be solved, their final GAP values after 6 hours ranged from 0.4\% to 6.2\%. The impact of SUB cuts is even higher as their removal results in 8 instances not being solved within 6 hours, with final GAP values ranging from 2.5\% to 27.6\%. All instances were solved for the other four cases (where we remove FCAV, CA, AVR and PT cuts). 

The average FGAP value is higher than the reinforced model FGAP for all cases. In fact, the model without subaisle cuts had overall higher FGAP values than even the original model itself. This can be explained by the separation of constraints (\ref{eq11}) for fractional solutions in FBC algorithm used for the original model.

For the cases where we remove FCAV only, CA only and PT only, the average times were actually lower than the average time for the reinforced model. This is mainly due to one particular instance, $\Delta = 20$ and $O = 20$, which the reinforced model solved in 11015.4 seconds, while the same instance was solved in 3940.5, 6386.9 and 6039.9 seconds by the model without FCAV, CA and PT cuts respectively. 


\begin{table}[!htbp]
\centering
{\scriptsize
\renewcommand{\tabcolsep}{1.05mm} 
\renewcommand{\arraystretch}{1.5} 
\begin{tabular}{|c|rrr|rrr|rrr|rrr|}
\hline
\multirow{2}{*}{Model} & \multicolumn{3}{c|}{$\Delta = 5, O \leq 20$} & \multicolumn{3}{c|}{$\Delta = 10, O \leq 20$} & \multicolumn{3}{c|}{$\Delta = 20, O \leq 15$} & \multicolumn{3}{c|}{Average}\\
 & \multicolumn{1}{c}{T(s)} & \multicolumn{1}{c}{GAP} & \multicolumn{1}{c|}{FGAP} & \multicolumn{1}{c}{T(s)} & \multicolumn{1}{c}{GAP} & \multicolumn{1}{c|}{FGAP} & \multicolumn{1}{c}{T(s)} & \multicolumn{1}{c}{GAP} & \multicolumn{1}{c|}{FGAP} & \multicolumn{1}{c}{T(s)} & \multicolumn{1}{c}{GAP} & \multicolumn{1}{c|}{FGAP} \\
\hline
\multicolumn{1}{|l|}{Reinforced model} 
   & 375.0  & 0.0  & 12.6 & 1296.7  & 0.0 & 17.4 & 867.6   & 0.0  & 24.5 & \textbf{845.8}   & \textbf{0.0} & \textbf{18.0} \\
\multicolumn{1}{|l|}{With A, SUB and AVR cuts}
   & 336.1  & 0.0  & 14.1 & 751.2   & 0.0 & 18.2 & 1719.3  & 0.0  & 25.4 & \textbf{913.1}   & \textbf{0.0} & \textbf{19.1} \\
\multicolumn{1}{|l|}{With A and SUB cuts}
   & 866.1  & 0.0  & 16.2 & 2633.3  & 0.0 & 19.3 & 3382.2  & 0.1  & 27.8 & \textbf{2262.8}  & \textbf{0.0} & \textbf{20.9} \\
\multicolumn{1}{|l|}{Original model}
   & 6756.7 & 3.5  & 20.8 & 12182.0 & 7.7 & 26.9 & 15312.2 & 13.7 & 37.6 & \textbf{11305.7} & \textbf{8.1} & \textbf{28.2} \\
\hline

\end{tabular}
}
\caption{Comparison of the reinforced model, original model and the reinforced model with selected valid inequalities. Only instances solved to optimality with the reinforced model are considered.}
\label{table8}
\end{table}

We also analysed the model without a selection of the least influential cuts. These results are shown in Table \ref{table8}. The second row in that table shows results for the model with A, SUB and AVR cuts only (so without FCAV, CA and PT cuts). The combined removal of these cuts actually causes the average time spent in solving each instance to be higher than the reinforced model. Also, the FGAP value, at 19.1\%, is higher than removing each individually. However, all instances were still solved within 6 hours. 

When we further remove AVR cuts (so the model with only A and SUB cuts), the average time and FGAP values, shown in the third row, are both considerably higher. Also, two instances could not be solved within 6 hours: $O = 20, \Delta = 10$ and $O = 15, \Delta = 20$. The results are also considerably worse when compared to the model with all but AVR cuts. These results can be partly explained as the addition of certain cuts may cause other cuts to be violated (as mentioned throughout Section~\ref{sec:ineq}), so the cuts cannot considered independently.

\subsection{Solving batching and routing separately}

From Table \ref{table3} it is clear that the largest instances are still challenging - after 6 hours, the optimality gap for $\Delta = 20, O = 30$ was 45.5\%. From a realistic standpoint, the reinforced model is not yet practical enough for use by large supermarkets which deal with many online orders a day. 

What happens in practice is that supermarkets assign orders to trolleys - perhaps via a heuristic, perhaps manually - and then compute routes individually - also heuristically or based simply on experience. However, even though optimal solutions to the joint problem are still not within reach for large instances, the reinforced model introduced in this paper can still be of benefit if the problems are treated separately. 

The  cuts introduced in Section~\ref{sec:ineq} tackle the routing part of the problem (as opposed to batching), and when the batching and routing problems are separated, \textbf{\emph{the proposed method is able to solve the routing problem to optimality very quickly}}. In this subsection we illustrate how we can take advantage of optimal routing. In particular we consider instances involving up to $O=5000$ orders.

\subsubsection{Small instances, larger warehouse: optimal routing}

Taking the largest instance, $\Delta = 20$ and $O = 30$, we know from Table \ref{table2} that 6 trolleys are sufficient to collect all orders. We generated 10 random assignments of orders to each of the 6 trolleys, making sure that capacity constraints are not violated. In each case, all 6 routes are solved separately and sequentially (parallelism was not explored). These results are summarised in Table \ref{table9}. The first row in the table shows that for the 10 random assignments, the average time to optimally solve for all 6 routes was only 5.4s, while in the worst case the time was around 7s. If there were more than 30 orders, the time necessary to solve all routes would grow linearly - more trolleys may be necessary, but we are still computing only single trolley routes.

\begin{table}[!htbp]
\centering
{\scriptsize
\renewcommand{\tabcolsep}{1.9mm} 
\renewcommand{\arraystretch}{1.5} 
\begin{tabular}{|cc|rrr|}
\hline
\multirow{2}{*}{Aisles} & \multirow{2}{*}{Cross-aisles} & \multicolumn{3}{c|}{Time in seconds} \\
 & & Mean & Max & Min \\
 \hline
 8  & 3 & 5.4  & 7.0  & 3.9  \\
 8  & 4 & 7.6  & 10.9 & 4.8  \\
 8  & 5 & 10.6 & 16.0 & 5.7  \\
12  & 3 & 26.9 & 49.1 & 11.5 \\
12  & 4 & 40.7 & 71.5 & 19.0 \\
12  & 5 & 64.0 & 90.4 & 25.8 \\
\hline
\end{tabular}
}
\caption{For the instance with $\Delta = 20, O = 30$, we randomly assign orders to trolleys (respecting trolley capacities) 10 different times. The number of trolleys required to collect all 30 orders is $6$, in each case all $6$ routes are solved independently, but sequentially (no parallel computing is used).}
\label{table9}
\end{table}

We also tested how fast routes are computed when the warehouse contains more aisles and more cross aisles. Having more aisles and cross-aisles means that the corresponding graph will have more artificial vertices and less products per subaisle. A higher proportion of vertices will have an outdegree of up to 4 and less vertices will have an outdegree of only 2, making the problem potentially harder. The remaining rows in Table \ref{table9} shows results for up to 5 cross-aisles as well as a warehouse with 12 aisles.

Average, minimum and maximum times are higher - for a warehouse with 12 aisles and 5 cross-aisles, 90s were needed in the worst case to solve all 6 routes, which is still a very competitive time in a real world setting. Given an assignment, defined perhaps heuristically, optimal routes can be computed quite quickly for large warehouses. 

\subsubsection{Large instances: adapting the heuristic to include optimal routing}

In Section~\ref{sec:heuristic} we described the time savings heuristic, used for warm-starting the exact algorithms. In that heuristic, routing is estimated by choosing the best out of four methods (s-shape, largest gap, combined and combined+). As the savings matrix is not recomputed, the algorithm requires $2O + O^2$ estimates for partial routes. Among these, $O$ estimates are required for computing individual trolley routes once final assignments are defined. 

Here we examine the benefit of replacing routing heuristics with optimal routing, for that we employ the IBC algorithm and solve the reinforced model. We compare three methods: (i) the original heuristic, as described in Section~\ref{sec:heuristic}, (ii) a heuristic with optimal routing for the final assignment, that is, we use estimates during the algorithm but, once all orders have been assigned to trolleys, we solve for each trolley individually to find its optimal route, and (iii) a heuristic with complete optimal routing, in which all $2O + O^2$ estimates are replaced by optimal routing.

Table \ref{table10} shows average solution values and the average elapsed times for $\Delta = \{5, 10, 20\}$. The last row is an aggregation of all instances considered. For methods (ii) and (iii), we added a third column with the percentage reduction in solution values as compared to the original heuristic.

\begin{table}[!htbp]
\centering
{\scriptsize
\renewcommand{\tabcolsep}{1mm} 
\renewcommand{\arraystretch}{1.4} 
\begin{tabular}{|c|rr|rrr|rrr|}
\hline
\multirow{2}{*}{$\Delta$} 
& \multicolumn{2}{c|}{(i)} 
& \multicolumn{3}{c|}{(ii)} 
& \multicolumn{3}{c|}{(iii)} \\

& \multicolumn{1}{c}{Value} & \multicolumn{1}{c|}{Time (s)} & \multicolumn{1}{c}{Value} & \multicolumn{1}{c}{Time (s)} & \multicolumn{1}{c|}{Red (\%)} & \multicolumn{1}{c}{Value} & \multicolumn{1}{c}{Time (s)} & \multicolumn{1}{c|}{Red (\%)} \\
\hline
5          &  749.2 & 0.0 &  721.9 & 0.5 & 3.7 &  683.0 &  6.1 & 8.7 \\
10         &  847.4 & 0.0 &  835.5 & 2.0 & 1.5 &  815.7 & 11.7 & 3.3 \\
20         & 1089.3 & 0.0 & 1066.9 & 1.8 & 2.1 & 1061.5 & 13.2 & 2.5 \\
Aggregated &  895.3 & 0.0 &  874.8 & 1.4 & 2.4 &  853.4 & 10.3 & 4.8 \\
\hline
\end{tabular}
}
\caption{Average solution values and elapsed times for three different heuristic configurations.}
\label{table10}
\end{table}

As methods (i) and (ii) differ only in the computation of routes once the final assignment is defined, the latter always produces better or equal results when compared to the former, on average (ii) produced solution values that were 2.4\% lower than (i). Method (iii), on the other hand, is not guaranteed to produce better or equal solutions than (i). However, solution values were on average 4.8\% lower. For the instances considered (up to $O = 30$), elapsed times are manageable even for (iii) (complete optimal routing). For example, the problem with $\Delta = 20, O = 30$ was solved in 39s using method (iii).

We also created very large instances in order to test the scalability of the algorithms. We created instances with $\Delta = 1$ and $O = \{50, 100, 200, 500, 1000, 2000, 5000\}$. That is, we gathered the first chronological $O$ orders of the Foodmart database into a single test instance. 
These large test instances are also publicly available.
The results for each individual instance are presented in Table \ref{table11}.

\begin{table}[!htbp]
\centering
{\scriptsize
\renewcommand{\tabcolsep}{1.2mm} 
\renewcommand{\arraystretch}{1.4} 
\begin{tabular}{|c|rr|rrr|rrr|}
\hline
\multirow{2}{*}{$O$} 
& \multicolumn{2}{c|}{(i)} 
& \multicolumn{3}{c|}{(ii)} 
& \multicolumn{3}{c|}{(iii)} \\

& \multicolumn{1}{c}{Value} & \multicolumn{1}{c|}{Time (s)} & \multicolumn{1}{c}{Value} & \multicolumn{1}{c}{Time (s)} & \multicolumn{1}{c|}{Red (\%)} & \multicolumn{1}{c}{Value} & \multicolumn{1}{c}{Time (s)} & \multicolumn{1}{c|}{Red (\%)} \\
\hline
50          &   2317.4 &   0.1 &   2231.2 &    3.2 & 3.7 &  2182.7 &    61.2 & 5.8 \\
100         &   4106.0 &   0.4 &   4012.1 &    5.5 & 2.3 &  3830.9 &   226.0 & 6.7 \\
200         &   6811.4 &   1.6 &   6653.7 &    6.5 & 2.3 &  6395.3 &   689.4 & 6.1 \\
500         &  14682.4 &  10.1 &  13899.5 &   16.1 & 5.3 & 13794.5 &  4486.6 & 6.0 \\
1000        &  27081.4 &  39.6 &  26160.7 &   53.7 & 3.4 &      -- &      -- & --  \\
2000        &  51689.0 & 153.1 &  49740.0 &  179.0 & 3.8 &      -- &      -- & --  \\
5000        & 122103.3 & 998.1 & 118269.7 & 1048.0 & 3.1 &      -- &      -- & --  \\
\hline
\end{tabular}
}
\caption{Solution values and elapsed times for large test instances and for three different heuristic configurations.}
\label{table11}
\end{table}

For $O > 500$, method (iii) is no longer a viable option due to very high computation times. However, up to $O = 500$, we observed an average reduction in solution values of 6.2\%. Method (ii), on the other hand, is a viable alternative as the increase in computation time is negligible and solution values are on average 3.4\% lower than (i).

\subsection{Investigating the tradeoff between waiting for more orders and total routing distance}
In practical situations as time passes new orders arrive. Potentially the more orders that are batched together to be simultaneously assigned to trolleys and routed the lower the total distance associated with routing all orders will be. It is possible to use the approach given in this paper to investigate this tradeoff between the number of orders batched together (which will involve waiting for new orders to appear) and total routing distance. 

To illustrate this we took two instances with $O=20$ for which,
from Table~\ref{table3}, we already know the optimal solution when all orders are batched together. We sequenced these orders so that we have an order of appearance (i.e.~we know the sequence in which orders appeared). We then applied the following algorithm: \\
\noindent Repeat until all orders have been assigned to a trolley:
\begin{itemize}[noitemsep,nolistsep]
\item  Take the first $K$ orders in the sequenced list of orders and solve the problem to optimality using the available trolleys. 
\item  The orders in the most heavily utilised trolley (the one involving the most baskets) will be regarded as assigned, so that these orders are removed from the sequenced list and the trolley removed from the problem (i.e.~we have decided to assign the orders in this trolley to it, so both the orders and the trolley can be removed from the problem). 
\end{itemize} 

\noindent At the end of this algorithm we will have a total routing distance associated with taking $K$ orders at a time; varying $K$ gives insight into the tradeoff between the number of orders batched together and total routing distance.

Table~\ref{table13} shows the results obtained for two instances with $O=20$ and $\Delta = 5,10$. The results for $K=20$ in that table  (taken from Table~\ref{table3}) show the total routing distance associated with batching all the orders together. As would be expected we have that (in general) as $K$ increases and we have more orders batched together to be simultaneously assigned to trolleys and routed computation time increases, but total routing distance decreases.  For example in that table for the problem with $\Delta = 5$ total routing distance decreases from 1295.5 when $K=5$ to 1014.5 when $K=10$ and to 870.4 when $K=20$.

\begin{table}[!htbp]
\centering
{\scriptsize
\renewcommand{\tabcolsep}{1.5mm} 
\renewcommand{\arraystretch}{1.4} 
\begin{tabular}{|rrr|rrr|}
\hline
\multicolumn{3}{|c|}{$\Delta = 5, O = 20$} & \multicolumn{3}{c|}{$\Delta = 10, O = 20$} \\
$K$ & UB & T(s) & $K$ & UB & T(s) \\
\hline
5        &     1295.5  &    11.3  &       5  &     1415.6  &    18.9  \\
6        &     1267.5  &    12.5  &       6  &     1311.6  &    18.6  \\
7        &     1012.6  &    17.0  &       7  &     1102.9  &    51.3  \\
8        &     1040.6  &    24.3  &       8  &     1082.5  &    79.9  \\
9        &     1020.5  &    51.7  &       9  &     1088.5  &   503.7  \\
10       &     1014.5  &    70.6  &      10  &     1076.6  &   790.2  \\
11       &     1024.9  &   183.0  &      11  &     1086.5  &   677.4  \\
12       &      988.9  &   396.9  &      12  &     1090.4  &  1596.0  \\
13       &      902.5  &   632.8  &      13  &     1016.8  &  4356.3  \\
14       &      936.5  &  2971.9  &      14  &     1032.4  &  4742.5  \\
15       &      940.5  &   656.9  &      15  &     1004.4  & 22986.2  \\
{\bf 20} & {\bf 870.4} & {\bf 2990.9} & {\bf 20} & {\bf 992.4} & {\bf 11015.4} \\
\hline
\end{tabular}
}
\caption{Routing distance tradeoff}
\label{table13}
\end{table}

We would make the point here that this investigation as to the tradeoff between waiting for more orders and total routing distance could be done in the same fashion as indicated here for much larger problems if we were to adopt an approach which uses heuristic batching  but optimal routing, such as seen in  the previous subsection.

Note here that a similar approach to that seen above can be used to investigate item placement, i.e.~where in the warehouse should items be stored, including in possible multiple locations if so desired. Simply use our approach with alternative placements (and example orders) and see which placement has the lowest total travel distance.

\subsection{No-reversal case}
\label{sec:noreversal}

A common feature to solutions of the problem as defined so far is trolley reversal in an aisle. By this we mean that a trolley might proceed (say north to south) down a subaisle, but then reverse direction and go back the subaisle to the cross-aisle vertex from which it came.  A special case of the problem as considered in this paper therefore is the no-reversal case, where a trolley is not allowed to reverse in a subaisle (however it can reverse at cross-aisle vertices). The motivation for this simplification is that solutions with no-reversal may be easier to implement in practice, as they may seem more logical and intuitive for humans.

In order to consider this special case we define a reduced graph with just one vertex in each subaisle, positioned at the midpoint of the subaisle, where that vertex could be used to pick all of the items in the subaisle. In term of the notation we have used this is equivalent to defining $R(a,c)=1~\forall a \in \mathbb{A}, \forall c \in \mathbb{C}$. In addition, since we now cannot reverse in a subaisle, we can set:
\begin{optprog}
& x_{t,v(a,c),S(a,c,1)} & \;=\; & x_{t,S(a,c,1),v(a,c+1)}  &\;\;\;\;\;\; \forall a \in \mathbb{A}, c \in \mathbb{C}\setminus[W_C] \label{jebnr1} \\
& x_{t,v(a,c+1),S(a,c,1)} & \;=\; & x_{t,S(a,c,1),v(a,c)}  &\;\;\;\;\;\; \forall a \in \mathbb{A}, c \in \mathbb{C}\setminus[W_C]  \label{jebnr2} 
\end{optprog}
These constraints ensure that if we use an arc to the vertex in the subaisle then we proceed onward to the next cross-aisle.

Clearly it might well be possible to produce a particularised formulation (e.g.~based upon arc routing) for this no-reversal case. However, since this special no-reversal case is not the primary focus of this paper, we have not explored this possibility. 

Computational results given in Table \ref{table12} indicate how our approach performs when adapted (as indicated above) for the no-reversal case. Here, the results for the reinforced model with reversal are as in Table \ref{table3}, but are repeated for convenience of comparison.

\begin{table}[!htbp]
\centering
{\scriptsize
\renewcommand{\tabcolsep}{0.75mm} 
\renewcommand{\arraystretch}{1.4} 
\begin{tabular}{|cc|rrrrrrr|rrrrrrr|rr|}
\hline
\multirow{2}{*}{$\Delta$} & \multirow{2}{*}{$O$} & \multicolumn{7}{c|}{Reinforced model} & \multicolumn{7}{c|}{Reinforced model with the special case of no-reversal} & \multicolumn{2}{c|}{Ratios (\%)} \\
 &  & \multicolumn{1}{c}{T(s)} & \multicolumn{1}{c}{UB} & \multicolumn{1}{c}{GAP} & \multicolumn{1}{c}{LB} & \multicolumn{1}{c}{FGAP} & \multicolumn{1}{c}{FLB} & \multicolumn{1}{c|}{NS} & \multicolumn{1}{c}{T(s)} & \multicolumn{1}{c}{UB} & \multicolumn{1}{c}{GAP} & \multicolumn{1}{c}{LB} & \multicolumn{1}{c}{FGAP} & \multicolumn{1}{c}{FLB} & \multicolumn{1}{c|}{NS} & \multicolumn{1}{c}{T(s)} & \multicolumn{1}{c|}{UB} \\
\hline
   5 &    5 &     3.9 &   348.6 & \multicolumn{1}{c}{--} & \multicolumn{1}{c}{--} &\textbf{  0.6} &    346.5 &    36 &\textbf{    0.5} &   384.6 & \multicolumn{1}{c}{--} & \multicolumn{1}{c}{--} &  2.1 &    376.5 &    77 &  12.8 &  10.3 \\
 &    6 &     2.9 &   364.8 & \multicolumn{1}{c}{--} & \multicolumn{1}{c}{--} &  0.0 &    364.8 &     1 &\textbf{    0.4} &   384.6 & \multicolumn{1}{c}{--} & \multicolumn{1}{c}{--} &  0.0 &    384.6 &     1  &  13.8 & 5.4  \\
 &    7 &     8.9 &   374.8 & \multicolumn{1}{c}{--} & \multicolumn{1}{c}{--} &  0.5 &    372.8 &   141 &\textbf{    0.9} &   384.6 & \multicolumn{1}{c}{--} & \multicolumn{1}{c}{--} &\textbf{  0.2} &    383.8 &     3 & 10.1 & 2.6 \\
 &    8 &    95.2 &   503.8 & \multicolumn{1}{c}{--} & \multicolumn{1}{c}{--} &  8.3 &    462.1 &   146 &\textbf{   10.9} &   543.7 & \multicolumn{1}{c}{--} & \multicolumn{1}{c}{--} &\textbf{  5.3} &    515.0 &    27 & 11.4 & 7.9 \\
 &    9 &   151.8 &   539.6 & \multicolumn{1}{c}{--} & \multicolumn{1}{c}{--} & 11.2 &    479.2 &   268 &\textbf{   25.5} &   603.2 & \multicolumn{1}{c}{--} & \multicolumn{1}{c}{--} &\textbf{  9.4} &    546.8 &   397 & 16.8 & 11.8 \\
 &   10 &   111.0 &   581.4 & \multicolumn{1}{c}{--} & \multicolumn{1}{c}{--} & 14.4 &    497.8 &    88 &\textbf{   21.1} &   611.4 & \multicolumn{1}{c}{--} & \multicolumn{1}{c}{--} &\textbf{  6.3} &    572.9 &   101 & 19.0 & 5.2 \\
 &   11 &    97.1 &   613.5 & \multicolumn{1}{c}{--} & \multicolumn{1}{c}{--} & 16.9 &    509.9 &  1009 &\textbf{   16.8} &   641.7 & \multicolumn{1}{c}{--} & \multicolumn{1}{c}{--} &\textbf{ 10.8} &    572.4 &    93 & 17.3 & 4.6 \\
 &   12 &   256.7 &   621.4 & \multicolumn{1}{c}{--} & \multicolumn{1}{c}{--} & 15.0 &    527.9 &  2261 &\textbf{    6.8} &   641.7 & \multicolumn{1}{c}{--} & \multicolumn{1}{c}{--} &\textbf{  9.6} &    579.9 &    85 & 2.6 & 3.3 \\
 &   13 &   168.9 &   623.4 & \multicolumn{1}{c}{--} & \multicolumn{1}{c}{--} & 15.3 &    528.2 &   395 &\textbf{   10.0} &   649.3 & \multicolumn{1}{c}{--} & \multicolumn{1}{c}{--} &\textbf{  9.2} &    589.3 &   114 & 5.9 & 4.2 \\
 &   14 &   263.8 &   639.3 & \multicolumn{1}{c}{--} & \multicolumn{1}{c}{--} & 17.8 &    525.5 &   995 &\textbf{   66.0} &   691.7 & \multicolumn{1}{c}{--} & \multicolumn{1}{c}{--} &\textbf{ 14.9} &    588.5 &  1221 & 25.0 & 8.2 \\
 &   15 &   348.9 &   653.4 & \multicolumn{1}{c}{--} & \multicolumn{1}{c}{--} & 20.0 &    522.4 &   975 &\textbf{   27.3} &   699.8 & \multicolumn{1}{c}{--} & \multicolumn{1}{c}{--} &\textbf{ 15.6} &    590.4 &  1710 & 7.8 & 7.1 \\
 &   20 &  2990.9 &   870.4 & \multicolumn{1}{c}{--} & \multicolumn{1}{c}{--} & 31.2 &    598.4 & 28931 &\textbf{  303.1} &   946.9 & \multicolumn{1}{c}{--} & \multicolumn{1}{c}{--} &\textbf{ 26.4} &    696.6 & 20694 & 10.1 & 8.8 \\
 &   25 & \multicolumn{1}{c}{--} &\textbf{  1127.4} &  17.8 &    927.0 &  43.2 &    640.8 & 123900 &\multicolumn{1}{c}{--} &  1155.1 & \textbf{  5.8} &   1088.6 & \textbf{ 34.2} &    760.1 & 327653 & \multicolumn{1}{c}{--} & \multicolumn{1}{c|}{--} \\
 &   30 & \multicolumn{1}{c}{--} &\textbf{  1221.9} &  26.8 &    894.3 &  47.0 &    647.5 & 70400 &\multicolumn{1}{c}{--} &  1295.0 & \textbf{ 16.4} &   1082.5 & \textbf{ 39.3} &    786.4 & 358800 & \multicolumn{1}{c}{--} & \multicolumn{1}{c|}{--} \\
\multicolumn{2}{|r|}{\textbf{Avg}} & \textbf{3407.1} &  & \textbf{3.2} &  & \textbf{17.2} &  \multicolumn{2}{r|}{\textbf{16396.1}} & \textbf{3120.7} &  & \textbf{1.6} &  & \textbf{13.1} &  \multicolumn{2}{r|}{\textbf{50784.0}} & \textbf{12.7} & \textbf{6.6}\\
\cline{1-18}
  10 &    5 &     1.5 &   371.1 & \multicolumn{1}{c}{--} & \multicolumn{1}{c}{--} &  0.0 &    371.1 &     1 &\textbf{    0.1} &   384.6 & \multicolumn{1}{c}{--} & \multicolumn{1}{c}{--} &  0.0 &    384.6 &     1 & 6.7 & 3.6 \\
 &    6 &     6.6 &   377.1 & \multicolumn{1}{c}{--} & \multicolumn{1}{c}{--} &  0.6 &    374.7 &   207 &\textbf{    0.4} &   384.6 & \multicolumn{1}{c}{--} & \multicolumn{1}{c}{--} &\textbf{  0.0} &    384.6 &     1 & 6.1 & 2.0 \\
 &    7 &   197.5 &   549.8 & \multicolumn{1}{c}{--} & \multicolumn{1}{c}{--} &  8.4 &    503.4 &    13 &\textbf{   16.2} &   613.1 & \multicolumn{1}{c}{--} & \multicolumn{1}{c}{--} &\textbf{  4.5} &    585.8 &   124 & 8.2 & 11.5 \\
 &    8 &   286.8 &   584.2 & \multicolumn{1}{c}{--} & \multicolumn{1}{c}{--} & 13.0 &    508.4 &   110 &\textbf{   22.0} &   681.2 & \multicolumn{1}{c}{--} & \multicolumn{1}{c}{--} &\textbf{  9.2} &    618.7 &  1497 & 7.7 & 16.6 \\
 &    9 &   356.1 &   637.4 & \multicolumn{1}{c}{--} & \multicolumn{1}{c}{--} & 17.1 &    528.6 &  8608 &\textbf{   17.1} &   681.2 & \multicolumn{1}{c}{--} & \multicolumn{1}{c}{--} &\textbf{  9.3} &    618.1 &   201 & 4.8 & 6.9 \\
 &   10 &   393.0 &   661.8 & \multicolumn{1}{c}{--} & \multicolumn{1}{c}{--} & 16.4 &    553.4 & 10439 &\textbf{   24.6} &   729.3 & \multicolumn{1}{c}{--} & \multicolumn{1}{c}{--} &\textbf{  9.3} &    661.7 &  3126 & 6.3 & 10.2 \\
 &   11 &   383.7 &   699.8 & \multicolumn{1}{c}{--} & \multicolumn{1}{c}{--} & 20.9 &    553.8 & 21044 &\textbf{   17.2} &   729.7 & \multicolumn{1}{c}{--} & \multicolumn{1}{c}{--} &\textbf{ 10.1} &    656.1 &   348 & 4.5 & 4.3 \\
 &   12 &   217.7 &   707.7 & \multicolumn{1}{c}{--} & \multicolumn{1}{c}{--} & 21.8 &    553.6 &  2481 &\textbf{   21.2} &   731.3 & \multicolumn{1}{c}{--} & \multicolumn{1}{c}{--} &\textbf{  9.1} &    664.7 &   515 & 9.7 & 3.3 \\
 &   13 &   385.4 &   725.7 & \multicolumn{1}{c}{--} & \multicolumn{1}{c}{--} & 22.8 &    560.2 & 10370 &\textbf{   98.5} &   769.2 & \multicolumn{1}{c}{--} & \multicolumn{1}{c}{--} &\textbf{ 13.1} &    668.7 &  6038 & 25.6 & 6.0 \\
 &   14 &   483.8 &   727.8 & \multicolumn{1}{c}{--} & \multicolumn{1}{c}{--} & 22.7 &    562.4 &  8597 &\textbf{  102.2} &   769.2 & \multicolumn{1}{c}{--} & \multicolumn{1}{c}{--} &\textbf{ 14.8} &    655.7 &  6755 & 21.1 & 5.7 \\
 &   15 &  1832.6 &   882.6 & \multicolumn{1}{c}{--} & \multicolumn{1}{c}{--} & 32.8 &    593.1 & 26398 &\textbf{  176.4} &   930.0 & \multicolumn{1}{c}{--} & \multicolumn{1}{c}{--} &\textbf{ 23.9} &    707.3 & 12528 & 9.6 & 5.4 \\
 &   20 & 11015.4 &   992.4 & \multicolumn{1}{c}{--} & \multicolumn{1}{c}{--} & 32.3 &    672.0 & 70496 &\textbf{  172.2} &  1028.3 & \multicolumn{1}{c}{--} & \multicolumn{1}{c}{--} &\textbf{ 21.5} &    807.5 & 15586 & 1.6 & 3.6 \\
 &   25 & \multicolumn{1}{c}{--} &\textbf{  1213.4} &  16.7 &   1011.3 &  41.7 &    707.8 & 133600 &\multicolumn{1}{c}{--} &  1264.7 & \textbf{  5.8} &   1190.7 & \textbf{ 32.8} &    849.5 & 276700 & \multicolumn{1}{c}{--}  & \multicolumn{1}{c|}{--} \\
 &   30 & \multicolumn{1}{c}{--} &\textbf{  1330.0} &  27.6 &    963.1 &  46.9 &    706.3 & 74100 &\multicolumn{1}{c}{--} &  1400.5 & \textbf{ 18.4} &   1142.8 & \textbf{ 40.2} &    836.9 & 318200 & \multicolumn{1}{c}{--} & \multicolumn{1}{c|}{--} \\
\multicolumn{2}{|r|}{\textbf{Avg}} & \textbf{4197.2} &  & \textbf{3.2} &  & \textbf{21.2} &  \multicolumn{2}{r|}{\textbf{26176.0}} & \textbf{3133.4} &  & \textbf{1.7} &  & \textbf{14.1} &  \multicolumn{2}{r|}{\textbf{45830.0}} & \textbf{9.3} & \textbf{6.6}\\
\cline{1-18}
  20 &    5 &   110.9 &   573.8 & \multicolumn{1}{c}{--} & \multicolumn{1}{c}{--} &  6.7 &    535.3 &   233 &\textbf{    8.5} &   623.3 & \multicolumn{1}{c}{--} & \multicolumn{1}{c}{--} &\textbf{  0.0} &    623.3 &     1 & 7.7 & 8.6 \\
 &    6 &   278.9 &   656.2 & \multicolumn{1}{c}{--} & \multicolumn{1}{c}{--} & 11.5 &    580.9 &  1164 &\textbf{   11.6} &   729.3 & \multicolumn{1}{c}{--} & \multicolumn{1}{c}{--} &\textbf{  6.9} &    679.3 &   734  & 4.2 & 11.1 \\
 &    7 &   263.2 &   689.8 & \multicolumn{1}{c}{--} & \multicolumn{1}{c}{--} & 17.2 &    571.4 &   579 &\textbf{   24.9} &   769.2 & \multicolumn{1}{c}{--} & \multicolumn{1}{c}{--} &\textbf{  9.2} &    698.1 &  3678 & 9.5 & 11.5 \\
 &    8 &   189.0 &   697.8 & \multicolumn{1}{c}{--} & \multicolumn{1}{c}{--} & 18.5 &    568.7 &  2886 &\textbf{   18.3} &   769.2 & \multicolumn{1}{c}{--} & \multicolumn{1}{c}{--} &\textbf{ 11.6} &    679.9 &  2795 & 9.7 & 10.2 \\
 &    9 &   305.5 &   727.7 & \multicolumn{1}{c}{--} & \multicolumn{1}{c}{--} & 20.6 &    578.0 &  3620 &\textbf{   49.4} &   769.2 & \multicolumn{1}{c}{--} & \multicolumn{1}{c}{--} &\textbf{ 10.4} &    689.3 &  1293 & 16.2 & 5.7 \\
 &   10 &   584.2 &   920.5 & \multicolumn{1}{c}{--} & \multicolumn{1}{c}{--} & 31.1 &    634.2 &  7198 &\textbf{  115.9} &   990.0 & \multicolumn{1}{c}{--} & \multicolumn{1}{c}{--} &\textbf{ 20.1} &    790.8 &  3435 & 19.8 & 7.6 \\
 &   11 &   703.8 &   980.5 & \multicolumn{1}{c}{--} & \multicolumn{1}{c}{--} & 31.4 &    672.9 & 10620 &\textbf{  132.4} &  1075.8 & \multicolumn{1}{c}{--} & \multicolumn{1}{c}{--} &\textbf{ 22.4} &    834.6 &  8361 & 18.8 & 9.7 \\
 &   12 &   830.1 &  1004.3 & \multicolumn{1}{c}{--} & \multicolumn{1}{c}{--} & 33.7 &    666.0 & 20221 &\textbf{  241.6} &  1113.8 & \multicolumn{1}{c}{--} & \multicolumn{1}{c}{--} &\textbf{ 23.9} &    847.9 & 29783 & 29.1 & 10.9 \\
 &   13 &  1177.8 &  1009.1 & \multicolumn{1}{c}{--} & \multicolumn{1}{c}{--} & 31.0 &    696.5 & 30072 &\textbf{  187.3} &  1113.8 & \multicolumn{1}{c}{--} & \multicolumn{1}{c}{--} &\textbf{ 23.5} &    851.7 & 28029 & 15.9 & 10.4 \\
 &   14 &  1306.0 &  1011.1 & \multicolumn{1}{c}{--} & \multicolumn{1}{c}{--} & 32.8 &    679.2 &  8097 &\textbf{  467.5} &  1113.9 & \multicolumn{1}{c}{--} & \multicolumn{1}{c}{--} &\textbf{ 28.4} &    797.6 & 23297 & 35.8 & 10.2 \\
 &   15 &  3793.7 &  1028.7 & \multicolumn{1}{c}{--} & \multicolumn{1}{c}{--} & 34.7 &    671.5 & 38243 &\textbf{  348.7} &  1115.9 & \multicolumn{1}{c}{--} & \multicolumn{1}{c}{--} &\textbf{ 28.4} &    798.7 & 24610 & 9.2 & 8.5 \\
 &   20 & \multicolumn{1}{c}{--} &\textbf{  1335.2} &  15.8 &   1124.6 &  45.4 &    729.0 & 67600 &\multicolumn{1}{c}{--} &  1453.0 & \textbf{  5.9} &   1367.0 & \textbf{ 37.4} &    910.0 & 341800 & \multicolumn{1}{c}{--} & \multicolumn{1}{c|}{--}\\
 &   25 & \multicolumn{1}{c}{--} &\textbf{  1694.9} &  31.2 &   1166.8 &  51.4 &    824.5 & 67200 &\multicolumn{1}{c}{--} &  1769.5 & \textbf{ 17.1} &   1466.5 & \textbf{ 40.2} &   1058.4 & 387359 & \multicolumn{1}{c}{--} & \multicolumn{1}{c|}{--}\\
 &   30 & \multicolumn{1}{c}{--} &\textbf{  1966.9} &  45.5 &   1071.9 &  56.6 &    853.4 & 33400 &\multicolumn{1}{c}{--} &  2041.8 & \textbf{ 29.5} &   1440.3 & \textbf{ 47.0} &   1082.6 & 360245 & \multicolumn{1}{c}{--} & \multicolumn{1}{c|}{--}\\
\multicolumn{2}{|r|}{\textbf{Avg}} & \textbf{5310.2} &  & \textbf{6.6} &  & \textbf{30.2} &  \multicolumn{2}{r|}{\textbf{20795.2}} & \textbf{4743.3} &  & \textbf{3.8} &  & \textbf{22.1} &  \multicolumn{2}{r|}{\textbf{86815.7}}  & \textbf{16.0} & \textbf{9.5}\\
\cline{1-18}
\hline
\multicolumn{2}{r}{\textbf{Avg all}} & \textbf{4304.8} &  & \textbf{4.3} &  & \textbf{22.9} &  \multicolumn{2}{r}{\textbf{21122.5}}
 & \textbf{3665.8} &  & \textbf{2.4} &  & \textbf{16.4} &  \multicolumn{2}{r}{\textbf{61143.2}} & \textbf{12.6} & \multicolumn{1}{r}{\textbf{7.5}}
\end{tabular}
}
\caption{Comparison of the reinforced model with and without the no-reversal special case. IBC was used in both cases.}
\label{table12}
\end{table}

Two aspects are of interest in this case: performance improvements and how much larger are solution values for the no-reversal case when compared to the reversal case. For analysing these aspects we have added two extra columns to Table \ref{table12}: T(s) ratio, defined as 100[no-reversal T(s)]/[reversal T(s)], shows how faster the 
no-reversal case is solved when compared to the reversal case. The UB ratio is defined as 100[no-reversal UB]/[reversal UB] - 100, which gives the percentage increase in solution value for the no-reversal model compared to the reversal model. Only instances solved to optimality are assigned values in these columns, the average ratios also take into account these instances only.

For instances solved to optimality, solving the no-reversal case is considerably faster than solving the reversal case, the time required for the former being on average only 12.6\% of the time required for the latter. Additionally, the optimal solutions of the model without reversal are on average only 7.5\% longer than the corresponding optimal solutions for the model with reversal. Rarely are optimal solutions without reversals more than 12\% worse than corresponding optimal solutions with reversals. The GAP and FGAP values also drop considerably for the no-reversal case.

The improvement was however still not enough to prove optimality of the largest instances. However, the special case of no-reversal may be worthy of further investigation since optimal solutions may be good enough for real-life problems, both in terms of distance travelled as well as being more intuitive and easier to operate.

\section{Application}
\label{sec:application}

Clearly in this paper we have adopted one particular type of order picking (where a trolley is pushed by a worker around aisles to collect the items needed for the orders assigned to the trolley). Our consideration of this problem was motivated by discussions with a major UK supermarket involved with home delivery of groceries. In other industries/contexts different picking methods (e.g. zone picking, bucket brigade picking) exist. However we believe that an important application of the approach given in this paper is to supermarket home deliveries, as discussed  further below.

Picking in the manner considered in this paper is widely used by supermarkets involved in home delivery (e.g. see~\citep{hays2005}). Indeed readers of this paper might themselves have seen supermarket staff pushing trolleys around stores and picking items destined for home delivery. The nature and economics of home delivery, where picking locations need to be relatively close to customers and where the capital investment associated with larger automated facilities can only be justified in areas of high population  density, mean that human pickers with trolleys (based either in publicly accessible or dark stores) will continue to be a large proportion of such operations, both for supermarkets in the UK and world-wide. This is especially true for supermarkets with a national presence (such as many in the UK) whose home delivery operations operate nationwide. 

With regard to the advantages of applying optimisation  \citep{bartholdi2016} 
state: \enquote{warehouses that are most likely to benefit from pick-path optimization are those
that have many items, most of which are slow-moving, and customer orders of moderate
size}.
Orders from customers to supermarkets for home delivery fall into this category. Moreover the supermarket business is a relatively low profit margin business and so pressure always exists to reduce costs appropriately. Obviously the use of optimisation to more effectively utilise  pickers  can aid in making home delivery operations more cost efficient. 

Clearly the difficulties in applying optimisation in the supermarket situation relate to communicating to the human trolley operator the precise route to be followed and ensuring that the route given is not overly complex.
 Overcoming these difficulties relies on a combination of technology and appropriate routing.  In terms of technology a simple handheld device (such as already used by supermarket pickers in order to know which items to pick) can be easily adapted to display instructions such as \enquote{now go to the end of aisle, turn left and back up the next aisle}. In terms of appropriate routing note that above we considered the no-reversal case, in which a trolley is not allowed to reverse in a subaisle. Our belief is that routes of this kind are simply easier for pickers to understand and follow as compared to routes which demand reversal in a subaisle.

With regard to the specific order characteristics used in this paper we applied our best judgement to enable us to produce orders with  characteristics as close to online grocery orders as achievable within the constraint of using a publicly available dataset (Foodmart). Using such a dataset means that we are able to make all the instances used in this paper publicly available for use by future workers.

\section{Conclusions and future work}
\label{sec:conclusion}

In this work we investigated the joint order batching and picker routing problem in storage areas. According to our literature review, very few exact methods have been proposed by any other authors, none for the case of warehouses composed of multiple blocks. We introduced a formulation and a significant number of valid inequalities (cuts) for the problem. Our results showed greatly improved performance. In fact, many potentially disconnected solutions are naturally prevented due to the addition of these new cuts. Instances involving up to 20 orders were solved to proven optimality when we jointly considered order batching and picker routing. Instances involving up to 5000 orders were considered where order batching was done heuristically, but picker routing done optimally.

In future we plan to study new algorithms for both JOBPRP and the special case of no-reversal, introduced in Section \ref{sec:noreversal}. We believe that column generation based methods, such as branch-and-price or branch-and-cut-and-price (BCP), may be viable alternatives. \cite{fukasawa2006} and \cite{baldacci2008} developed successful BCPs for the closely related Capacitated VRP and \cite{letchford2009} proposed column generation methods that exploit sparsity for the also closely related Capacitated Arc Routing Problem (CARP). If non-elementary routes are allowed (routes where edges/vertices may be revisited), then the authors showed that exact pricing can be computed polynomially. In our opinion, a BCP approach for the JOBPRP should exploit its sparse nature and perhaps an effective pricing algorithm can also be obtained. Moreover, we also believe that a CARP based formulation may be beneficial for the special case of no-reversal, in which a similar method to that in \citep{letchford2009} may be achievable.

We also plan to conduct a theoretical study of the strength of each valid inequality and to potentially characterize new valid inequalities that could be added to further reinforce linear relaxation bounds. Finally, we intend to study assignment heuristics where optimal routing replaces (at least partially) approximate route estimates.


\bibliographystyle{plainnat}
\bibliography{tesco}

\end{document}